\documentclass[a4paper]{cas-sc}

\usepackage[authoryear,longnamesfirst]{natbib}
\usepackage{amsmath, amsfonts, amssymb, amsthm, latexsym}
\usepackage{mathtools} 

\usepackage{graphicx}
\usepackage[all,cmtip]{xy} 
\usepackage{tikz-cd}
\usepackage{subcaption}
\usepackage{graphicx}
\usepackage{subcaption}
\usepackage{footmisc}

\DeclareMathOperator{\re}{Re}
\DeclareMathOperator{\im}{Im}

\newcommand{\R}{\mathbb{R}}
\newcommand{\C}{\mathbb{C}}

\renewcommand{\H}{\mathbb{H}}

\newcommand{\Sp}{\mathrm{Sp}}

\newcommand{\bild}{\operatorname{bild}}

\def\tsc#1{\csdef{#1}{\textsc{\lowercase{#1}}\xspace}}
\tsc{WGM}
\tsc{QE}

\newtheorem{theorem}{Theorem}
\newtheorem{lemma}[theorem]{Lemma}
\newtheorem{proposition}[theorem]{Proposition}
\newtheorem{example}[theorem]{Example}

\newtheorem{definition}[theorem]{Definition}
\newdefinition{remark}{Remark}
\newproof{pf}{Proof}
\newproof{pot}{Proof of Theorem \ref{thm}}

\begin{document}
\let\WriteBookmarks\relax
\def\floatpagepagefraction{1}
\def\textpagefraction{.001}

\shorttitle{}    

\shortauthors{J. Tian and T.-Y. Tam}  

\title [mode = title]{Revisiting the Upper Bild of Quaternionic Numerical Ranges}  



%

\author[1]{Jie Tian}[orcid=0009-0002-8793-2304]

\cormark[1]

\ead{jtian@unr.edu}

\ead[url]{https://jietinatian0803.github.io}

\credit{Conceptualization, Methodology, Formal analysis, Validation, Software, Visualization, Writing -- original draft, Writing -- review \& editing}

            
 \affiliation[1]{organization={University of Nevada, Reno},
 	addressline={1664 N. Virginia Street}, 
 	city={Reno},
 	postcode={89557}, 
 	state={Nevada},
 	country={USA}}

\author[1]{Tin-Yau Tam}

\ead{ttam@unr.edu}

\ead[url]{https://sites.google.com/view/tinyautam/homepage}

\credit{Conceptualization, Supervision, Writing -- review \& editing}

\cortext[1]{Corresponding author}



\begin{abstract}
	The upper bild of the quaternionic numerical range provides a convex representative of the generally nonconvex quaternionic numerical range. Although its convexity is known, the geometry of the upper bild is difficult to describe explicitly, even for \(2\times2\) quaternionic matrices. 
	
	In this paper, we revisit the upper bild of quaternionic numerical ranges from both algebraic and geometric viewpoints. We first point out gaps in Kumar's proof. We then investigate several special classes of \(2\times2\)
	quaternionic matrices and describe their shapes and boundary
	structures in the closed upper complex half-plane directly in
	terms of the standard eigenvalues and the off-diagonal entry
	of the Schur form.
    
    For the general \(2\times2\) case, we use a different canonical form.
    Following So and Thompson, we prove the concavity of the upper-boundary
    function through the associated complex representation.
    For the lower boundary, we obtain a distance-to-a-disk representation
    and use convex optimization to prove its convexity.
    We also treat the degenerate cases.
    Together with the connectedness of the vertical sections, these
    results give an alternative proof of the convexity of the upper
    bild for arbitrary quaternionic matrices.
	
\end{abstract}




\maketitle

\section{Introduction}\label{Section: Introduction}

The numerical range is a classical object in matrix analysis and operator theory that encodes spectral and geometric information about a matrix.
For a complex matrix, the celebrated Toeplitz--Hausdorff theorem asserts that the numerical range is always a convex subset of the complex plane. In contrast, the numerical range of a quaternionic matrix need not be convex. Nevertheless, the upper bild of the quaternionic numerical range is always convex \cite{ST1}.

Let $q=a+bi+cj+dk \in\H$ and $A=(q)\in M_{1}(\H)$. The quaternionic numerical range of $A$ is 
\begin{equation*}
	W_{\H}(A)=\{ x^{*} A x: x\in \Sp(1)\} =\{ p \in \H : \re(p)=\re(q), |\im (p)|= |\im (q)|\}.
\end{equation*}
If \(q\notin\R\), then \(W_{\H}(A)\) is a \(2\)-sphere in \(\H\) centered at \(a\) with radius
\(
|\im(q)|=\sqrt{b^2+c^2+d^2},
\)
and hence is not convex. If \(q\in\R\), then \(W_{\H}(A)=\{q\}\). However, 
\(
\bild^+(A)
=
\left\{a+i\sqrt{b^2+c^2+d^2}\right\}
\) is a singleton and hence convex.

The quaternionic numerical range is invariant under conjugation by unit quaternions and is therefore a union of quaternionic similarity classes. The upper bild selects the unique representative of each class in the closed upper complex half-plane. In this way, it removes the conjugation symmetry and gives a planar representation of the quaternionic numerical range. This planar representation makes geometric properties such as convexity easier to study.

Let $\H$ denote quaternions, $M_n(\H)$ be the set of $n\times n$ matrices with quaternionic entries and
\(
\Sp(1)=\{\,q\in  \H:\ |q|=1\,\}
\) denote the group of unit quaternions,
where $|q|$ is the standard quaternion norm.  
This group may be identified with the unit $3$--sphere in $\mathbb R^4$. Two quaternions $p,q$ are similar if there is a $u\in \Sp(1)$ such that $u^{*}pu=q$. We write $p\sim q$. The equivalence class of $q$, denoted by $[q]_{\sim}$, is given by
\[
[q]_{\sim} = \{ p\in \H: p\sim q \}. 
\]
It is easy to verify that $\re(q)=\re (p)$ and $ |\im(q)|= |\im(p)|$ iff $p\sim q$. Geometrically, the equivalence class of $q$ is a $2$--sphere in $\H$ centered at $\re(q)$ with radius $|\im(q)|$ (degenerating to a point when $\im(q)=0$).
For $x\in\H^n$,
we write $x^*$ for the conjugate transpose of $x$ and $\|x\|=(x^*x)^{1/2}$ for the
standard quaternionic norm.  The unit sphere in $\H^n$ is denoted by
\[
\mathbb S^{n-1}_{\H}=\{x\in\H^n:\ \|x\|=1\}.
\]
For a given matrix $A\in M_n(\H)$, its quaternionic numerical range is
defined by
\[
W_{\H}(A)
=
\{\,x^*Ax:\ x\in\mathbb S^{n-1}_{\H}\,\}
\subset \H.
\]
The quaternionic numerical range $W_{\H}(A)$ is invariant under
conjugation by unit quaternions, that is, for every $u\in\Sp(1)$ and every
$x\in\mathbb S^{n-1}_{\H}$, if $x^*Ax=q \in W_{\H}(A)$, then
\[
u^{*}qu=u^{*}(x^*Ax)u=(xu)^*A(xu)\in  W_{\H}(A). 
\]

\begin{definition}
	Define the upper bild map
	\[
	\pi: \H \longrightarrow \C^+,
	\qquad
	\pi(q) = \re(q) + i\,|\im(q)|,
	\]
	where $\C^+ = \{ z\in\C : \im (z) \ge 0 \}$.
\end{definition}

\begin{definition}[Upper Bild]
	The upper bild of the quaternionic numerical range of $A\in M_{n}(\H)$ is
	\[
	\bild^{+}(A)=W_{\H}(A) \cap \C^{+}=\pi\bigl(W_{\H}(A) \bigr).
	\]
\end{definition}

Note that the quaternionic numerical range of $A$ can be expressed as
\[
W_{\H}(A)=\bigcup_{q\in \bild^{+}(A)} [q]_{\sim}.
\]

We next record some elementary properties of the quaternionic numerical range.

\begin{proposition}[\cite{K19}]

	Let $A\in M_n(\mathbb{H})$.  The following properties hold true:
	\begin{enumerate}[a)]
		\item
		\(
		W_{\H}(\alpha I+\beta A)
		=
		\alpha+\beta W_{\H}(A),
		\quad \text{for every } \alpha,\beta\in\mathbb{R}.
		\)
		
		\item
		If $B\in M_n(\mathbb{H})$, then
		\(
		W_{\mathbb{H}}(A+B)
		\subseteq
		W_{\mathbb{H}}(A)+W_{\mathbb{H}}(B).
		\)
		
		\item
		\(
		W_{\mathbb{H}}(U^{*}AU)
		=
		W_{\mathbb{H}}(A),
		\quad
		\text{for every unitary } U\in M_n(\mathbb{H}).
		\)
		
		\item
		\(
		W_{\mathbb{H}}(A^{*})
		=
		W_{\mathbb{H}}(A).
		\)

	\end{enumerate}
	
\end{proposition}

\begin{proposition}
		The upper bild has the following properties:
	\begin{enumerate}
		\item  	In general, there exist \(A,B\in M_n(\H)\) such that
		$
		\bild^{+}(A+B)
		\not\subseteq
		\bild^{+}(A)+\bild^{+}(B).
		$
		
		\item 
		For every \(A\in M_n(\H)\), \(\alpha\in\R\), and \(\beta\geq0\),
		\(
		\bild^{+}(\alpha I+\beta A)
		=
		\alpha+\beta\bild^{+}(A).
		\)
	\end{enumerate}
	
\end{proposition}

 \begin{example}
	Let
	\[
	A=\begin{pmatrix}
		1+2i&0\\
		0&3+i
	\end{pmatrix},
	\qquad
	B=\begin{pmatrix}0&2j\\0&0\end{pmatrix}.
	\]
	Taking \(x=(3,\,3j-k)^T/\sqrt{19}\), one readily verifies that
	\[
	z=x^*(A+B)x=\frac{21+2i}{19}
	\in\operatorname{bild}^{+}(A+B).
	\]
	On the other hand, every \(r+yi\in
	\operatorname{bild}^{+}(A)+\operatorname{bild}^{+}(B)\)
	satisfies \(r+y\geq 4/3\), while
	\(\operatorname{Re}z+\operatorname{Im}z=23/19<4/3\).
	Thus
	\[
	\operatorname{bild}^{+}(A+B)
	\not\subseteq
	\operatorname{bild}^{+}(A)+\operatorname{bild}^{+}(B).
	\]
\end{example}
Although the quaternionic numerical range enjoys several algebraic properties such as Minkowski-type inclusion, unitary invariance, and affine equivariance, the upper bild should be regarded as a nonlinear orbit representative rather than a linear projection.  In general, taking the upper bild of the numerical range of a quaternionic matrix does not commute with algebraic operations on matrices.

 \cite {ST1} established the convexity of $\bild^+(A)$. Their proof first reduces the general
\(n\times n\) case to the \(2\times2\) case. They then divide the
\(2\times2\) case into several special cases and one generic case, and
study them separately using explicit parametrizations and detailed
calculations. For the upper boundary, So and Thompson used a result implicit in \cite{K51},
which states that the upper boundary of the bild coincides with the
upper boundary of the complex numerical range of its complex
representation. Since the latter is convex by the
Toeplitz--Hausdorff theorem, its upper boundary function is concave. In the generic \(2\times2\) case, the lower boundary is divided into
three parts: a left nonreal branch, a line segment on the real axis,
and a right nonreal branch. They proved that the derivative of the
lower boundary function is monotonically increasing on each of the two
nonreal branches. Since the lower boundary function is nonnegative and
vanishes on the middle real-axis interval, the three parts fit together
to form a convex lower boundary function. Together with the concavity of the upper boundary function and the connectedness of the vertical sections, this proves the convexity of the upper bild in the generic \(2\times2\) case. The special \(2\times2\) cases are treated separately. Combining these
special cases with the generic case proves the convexity of the upper
bild for every \(2\times2\) quaternionic matrix. Finally, by the
reduction from the general \(n\times n\) case to the \(2\times2\) case,
the convexity of the upper bild follows for every quaternionic matrix.

\cite{K19} later claimed to give a new proof of the convexity of the upper bild. However, there are several gaps in the proof. We explain them in Section~\ref{Gaps in Kumar's paper}.

Recently, \cite{carvalho2026convexity} gave an independent
proof of the convexity of the upper bild. His proof introduces real
coordinates on a sphere and extends two representations of the upper
bild to the closed unit ball without changing their images.
The two radius functions agree on the sphere, while their extensions
to the ball are convex and concave, respectively. Minimizing the
convex function and maximizing the concave function over affine
fibres with a fixed real part yield a convex lower boundary and a
concave upper boundary. This gives a concise and conceptual proof
of the convexity of the upper bild.

In this paper, we study the geometry of the upper bild for several
special classes of \(2\times2\) quaternionic matrices and give an
alternative proof of its convexity for arbitrary quaternionic matrices. We also begin by reducing the general \(n\times n\)
problem to the \(2\times2\) case. We also use a quotient-space description of quaternionic
unit vectors to remove the redundancy caused by right multiplication
by unit quaternions. To study the geometry of the special classes, we use quaternionic
unitary similarity to put a \(2\times2\) matrix in upper triangular form and reduce the number of parameters. The upper triangular form allows us to describe the geometry of these special classes in terms of the matrix parameters. As a by-product, we determine when the quaternionic numerical ranges of matrices in these classes are convex.
For the general \(2\times2\) case, we use a different canonical form
to prove the convexity of the upper bild. Following So and Thompson,
we relate the upper boundary to the complex numerical range of the
associated complex representation and give a detailed proof of its
concavity. For the lower boundary, we use convex duality to obtain a
distance-to-a-disk representation in a fixed positive-definite
quadratic metric and prove that the lower-boundary function is convex.
Degenerate cases are treated directly or by approximation.
Together with the connectedness of the vertical sections, these
results establish the convexity of the upper bild for every
\(2\times2\) quaternionic matrix, and hence in every dimension.
This distance-based argument differs from Carvalho's
sphere-to-ball construction \cite{carvalho2026convexity}.

The remainder of this paper is organized as follows.
Section~\ref{Gaps in Kumar's paper} identifies gaps in Kumar's proof.
Section~\ref{Section:Reduction} reduces the general problem to the
\(2\times2\) case and introduces the quotient-space framework used
throughout the paper. Section~\ref{Section:Specialshape} studies several special classes of
\(2\times2\) quaternionic matrices.
Section~\ref{Section:Generic} treats the general \(2\times2\) case and
completes the proof of the convexity of the upper bild.
Finally, Section~\ref{Section:Conclusion} summarizes the main results
and discusses possible directions for future research.

\section{\textbf{Gaps in Kumar's paper}}\label{Gaps in Kumar's paper}

In this section, we use the notation \(W_{\H}^{+}(A)\) from \cite{K19}; in our notation, this set is \(\bild^{+}(A)\). Also, here \(\Omega_S\) denotes the circularization of \(S\subseteq\C\), as in \cite{K19}.

In Case~(2) of \cite[Lemma~2.14]{K19} in Kumar's paper,  the proof uses the implication
\[
z \in \mathbb{C}^+ \setminus \operatorname{conv}\{z_1,z_2,v\}
\quad \Longrightarrow \quad
z \notin \operatorname{conv}\bigl(\Omega_{\operatorname{conv}\{z_1,z_2,v\}}\bigr),
\]
which does not hold in general. Indeed, if $\operatorname{conv}\{z_1,z_2,v\} \subset \mathbb{C}^+$, 
then $\operatorname{conv}\{\overline{z_1},\overline{z_2},v\}\subset  \Omega_{\operatorname{conv}\{z_1,z_2,v\}}. $
As illustrated by the figure for Case~(2) on p.~108, there exist points
$z \in \mathbb{C}^+$ lying arbitrarily close to the line segment joining
$z_1$ and $v$ such that
\[
z \notin \operatorname{conv}\{z_1,z_2,v\},
\quad \text{but} \quad
z \in \operatorname{conv}\bigl(\Omega_{\operatorname{conv}\{z_1,z_2,v\}}\bigr).
\]
Therefore, the conclusion
\(
z \notin W_{\mathbb{H}}^+(A)
\)
cannot be deduced from the above exclusion. More specifically,
$$\operatorname{conv} ( \operatorname{conv} \{ z_{1}, z_{2}, v\}  \cup \operatorname{conv} \{ \overline{z_1},\overline{z_2},v\}) \setminus \operatorname{conv} \{ z_{1}, z_{2}, v\} $$
is nonempty, where $ \operatorname{conv} \bigl( \Omega_{\operatorname{conv}\{z_1,z_2,v\}} \bigr) \supset \operatorname{conv} ( \operatorname{conv} \{ z_{1}, z_{2}, v\}  \cup \operatorname{conv} \{ \overline{z_1},\overline{z_2},v\})$.

 In Case~(4) of \cite[Lemma~2.14]{K19}, the proof claims that the upper section \(W_\H^+(A)\) is a closed subset of a convex set \(\Gamma\), and then concludes that \(W_\H^+(A)\) is convex. This inference is incorrect. A closed subset of a convex set need not itself be convex. For example, \(\{1,2\}\) is a closed subset of the convex set
 \([1,2]\), but it is not convex, since
 \[
 \frac{1+2}{2}=\frac{3}{2}\notin\{1,2\}.
 \]
 
 	In Subcase~(4.3) of \cite[Lemma~2.14]{K19}, Kumar considers
 	\[
 	A=\begin{pmatrix}
 		3+4i & 1-j\\
 		0 & -2+5i
 	\end{pmatrix},
 	\]
 	that is, $z_1=3+4i$, $z_2=-2+5i$ and $p=1-j$, and asserts that
 	\[
 	\frac{|p|}{2}\,i=\frac{i}{\sqrt2}\notin W_{\H}^{+}(A).
 	\]
 	This assertion is incorrect. Take
 	\[
 	x=\frac{1}{\sqrt2}\begin{pmatrix}1\\-j\end{pmatrix}\in\mathbb S_{\H}^{1},
 	\qquad
 	x^*=\frac{1}{\sqrt2}\begin{pmatrix}1 & j\end{pmatrix}.
 	\]
 	Since
 	\[
 	(1-j)(-j)=-1-j,
 	\qquad
 	j(-2+5i)(-j)=-2-5i,
 	\]
 	we obtain
 	\[
 	x^*Ax
 	=\frac12\Bigl[(3+4i)+(1-j)(-j)+j(-2+5i)(-j)\Bigr]
 	=\frac12\Bigl[(3+4i)+(-1-j)+(-2-5i)\Bigr]
 	=-\frac12\,i-\frac12\,j .
 	\]
 	Hence $\re(x^*Ax)=0$ and $|\im(x^*Ax)|=\frac{1}{\sqrt2}$, so that \(\pi(x^*Ax)=\frac{i}{\sqrt2}\in\bild^{+}(A)=W_{\H}^{+}(A).\)

	In Case~(4) of \cite[Lemma~2.14]{K19}, for
	\[
	A=\begin{pmatrix} a_1+b_1i & p\\ 0 & a_2+b_2i\end{pmatrix},
	\qquad
	p=p_0+p_1i+p_2j+p_3k\neq0,
	\]
	Kumar considers the anti-Hermitian matrix
	\[
	B=\begin{pmatrix} b_1i & \frac p2\\[1mm] -\frac{\overline p}{2} & b_2i\end{pmatrix}
	\]
	and asserts that
	\[
	\det(\chi_B)=(b_1b_2)^2+\tfrac12\,b_1b_2|p|^2+\tfrac1{16}|p|^4>0,
	\]
	so that $B$ is always invertible. This formula is incorrect. A direct
	computation gives
	\[
	16\det(\chi_B)
	=16b_1^2b_2^2-8b_1b_2\,(p_0^2+p_1^2)+8b_1b_2\,(p_2^2+p_3^2)+|p|^4 .
	\]
	Since \(b_1,b_2\geq0\) for the chosen standard eigenvalues, set
	\(\mu=\frac14(p_0^2+p_1^2)\) and \(\nu=\frac14(p_2^2+p_3^2)\). Then
	\[
	\det(\chi_B)=(b_1b_2-\mu)^2+\nu\,(2b_1b_2+2\mu+\nu)\;\geq\;0,
	\]
	with equality if and only if $p_2=p_3=0$ and $|p|^2=4b_1b_2$.
	Hence $B$ can be singular. For instance, let
	\[
	A=\begin{pmatrix} i & 2\\ 0 & 1+i\end{pmatrix},
	\qquad\text{so that}\qquad
	B=\begin{pmatrix} i & 1\\ -1 & i\end{pmatrix}.
	\]
	Then
	\[
	B\begin{pmatrix}1\\-i\end{pmatrix}
	=\begin{pmatrix} i-i\\ -1+1\end{pmatrix}
	=\begin{pmatrix}0\\0\end{pmatrix},
	\]
	so $B$ is singular, whereas Kumar's formula would give $\det(\chi_B)=4>0$.

\begin{remark}
	We note that the degenerate value $|p|^2=4b_1b_2$ with
	$p\in\C$ is exactly the threshold $m^2=4b_1b_2$ found in
	Section~\ref{m>0}, at which the associated ellipse is tangent to the real axis.
\end{remark}

\section{\textbf{Reduction}}\label{Section:Reduction}

Like \cite{ST1} and \cite{K19}, we use the following
result to reduce the convexity problem from arbitrary dimension to the
\(2\times2\) case.

\begin{lemma}[\cite{K19},\cite{ST1}]\label{2x2}
	Let $A\in M_n(\H)$ be a quaternionic matrix, where \(n\geq2\).
	Assume that for every $2\times2$ quaternionic matrix $C$, the set $\bild^+(C)$ is convex.
	Then $\bild^+(A)$ is convex. 
\end{lemma}

Let $A\in M_{n}(\H)$. For any quaternionic unitary matrix $U\in M_{n}(\H)$, and for any unit vector \(x\in\mathbb{H}^{n}\), one has that
\[
x^{*}(U^{*}AU)x=(Ux)^{*}A(Ux)=y^{*}Ay,
\]
where \(y:=Ux\) is also a unit vector. Hence, it suffices to consider simpler cases after applying a suitable unitary similarity transformation.

\begin{lemma}[\cite{Z97}]\label{uptri}
	Let \(A\in M_n(\H)\). Then there exists a quaternionic unitary matrix
	\(U\in M_n(\H)\) such that \(U^*AU\) is upper triangular with diagonal
	entries
	\(
	a_\ell+b_\ell i,~\ell=1,\ldots,n,
	\)
	where \(a_\ell,b_\ell\in\R\) and \(b_\ell\geq0\). These diagonal entries
	are the standard eigenvalues of \(A\).
\end{lemma}

For the geometric analysis in Section~\ref{Section:Specialshape},
we use Lemma~\ref{uptri} to put the \(2\times2\) matrix in upper
triangular form with its standard eigenvalues on the diagonal.

We next remove the redundancy arising from right multiplication of unit vectors by unit quaternions. Recall that 
\[
\mathbb{HP}^1
:=
(\mathbb{H}^2\setminus\{0\})/\dot{\sim}, 
\]
where
\[
(x_1,x_2)^{T}\dot{\sim} (y_1,y_2)^{T}
\quad\Longleftrightarrow\quad
\exists\, q\in\mathbb H\setminus \{ 0\}
\text{ such that }
(y_1,y_2)^{T}=(x_1q,x_2q)^{T}.
\]
Equivalently, $\mathbb{HP}^1$ is the set of all right quaternionic
one-dimensional subspaces of $\mathbb{H}^2$. For unit vectors in $\mathbb H^2$, we have the following lemma. 
\begin{lemma}\label{xp}
	Let $x=(x_1,x_2)^T\in\mathbb H^2$ satisfy $\|x\|=1$.
	Then there exist $\theta\in[0,\tfrac{\pi}{2}]$ and $u\in \Sp(1)$ such that
	\[
	x \dot{\sim}
	\begin{pmatrix}
		\cos\theta\\
		\sin\theta\,u
	\end{pmatrix},
	\]
	where $\dot{\sim}$ denotes equivalence up to right multiplication by a unit quaternion.
\end{lemma}

\begin{proof}
	Let $x=(x_1,x_2)^T\in\mathbb H^2$ with $\|x\|^2=|x_1|^2+|x_2|^2=1$.
	We define $\theta \in [0, \frac{\pi}{2}]$ such that
	\[
	\cos\theta := |x_1|, \quad \sin\theta := |x_2|.
	\]
	To find the representative element, consider the polar decomposition of $x_1$. Let $v \in \Sp(1)$ be a unit quaternion such that $x_1 = |x_1|v$. Specifically, if $x_1 \neq 0$, $v = x_1/|x_1|$; if $x_1 = 0$, we can choose $v=1$.

	Now, right-multiply the vector $x$ by $v^*$:
	\[
	x v^* = \begin{pmatrix} x_1 v^* \\ x_2 v^* \end{pmatrix} 
	= \begin{pmatrix} (|x_1|v) v^* \\ x_2 v^* \end{pmatrix} 
	= \begin{pmatrix} |x_1| \\ x_2 v^* \end{pmatrix}
	= \begin{pmatrix} \cos\theta \\ x_2 v^* \end{pmatrix}.
	\]
	If $\sin\theta = 0$ (i.e., $x_2 = 0$), the second component is $0$, and the claim holds with any $u \in \Sp(1)$. 
	If $\sin\theta \neq 0$, we define:
	\[
	u := \frac{x_2 v^*}{|x_2 v^*|} = \frac{x_2 v^*}{\sin\theta} \in \Sp(1).
	\]
	It follows that $x_2 v^* = \sin\theta \, u$. Substituting this back, we obtain that
	\[
	x \dot{\sim} x v^* = \begin{pmatrix} \cos\theta \\ \sin\theta \, u \end{pmatrix}.
	\]
	Since $v^* \in \Sp(1)$, the vector $x v^*$ represents the same point in $\mathbb{HP}^{1}$ as $x$.
\end{proof}

We now describe the relationship between $\mathbb{S}_{\mathbb{H}}^{1}=\{ x| x\in \mathbb{H}^{2}, || x||=1\}$, $\mathbb{S}_{\mathbb{HP}}^1
:=
\mathbb S_{\H}^{1}/\Sp(1)$, $W_{\mathbb{H}}(A)$ and $\bild^+(A)$. Denote the induced map
\[
\rho_A:\mathbb S_{\mathbb{HP}}^{1}\to W_{\mathbb H}(A)/\!\sim,\qquad
[x]_{\dot{\sim}}\mapsto [x^*Ax]_{\sim}.
\]

A more precise description is the following:

\[
\mathbb S_{\H}^1
\;\xrightarrow{\ /\,\Sp(1)\ }\;
\mathbb S_{\mathbb{HP}}^{1}
\;\xrightarrow{\ \rho_A\ }\;
W_{\H}(A)/\!\sim
\;\cong\;
\bild^+(A).
\]

For a projective point $x$ in $\mathbb S_{\mathbb{HP}}^{1}$, the quantity \(x^*Ax\) depends on the chosen representative \(x\), but different representatives produce similar quaternions, that is, 
\[
(xu)^*A(xu)=u^*(x^*Ax)u,\qquad u\in \Sp(1).
\]
Therefore the induced map is well defined. Finally, each quaternionic similarity class meets the upper half-plane \(\mathbb C^+\) in exactly one point, so that
\[
W_{\mathbb H}(A)/\!\sim \;\cong\; \bild^+(A).
\]
Thus \(\bild^+(A)\) records the projective numerical range after quotienting out quaternionic similarity.

Moreover, in the Rayleigh quotient \(x^{*}Ax\), it suffices to choose one unit representative from each class in \(\mathbb S_{\mathbb{HP}}^{1}\).

\section{\textbf{The Shape of the Upper Bild for Special Classes of $2\times2$ Quaternionic Matrices}}\label{Section:Specialshape}

Let \(A\in M_2(\H)\) be upper triangular, and choose a unit representative
\(x\in\mathbb S_{\H}^{1}\) of a class in \(\mathbb S_{\mathbb{HP}}^{1}\), with
\[
x=
\begin{pmatrix}
	\cos\theta \\
	\sin\theta\, u
\end{pmatrix},~A=\begin{pmatrix}
c_{1}&q\\0&c_{2}
\end{pmatrix},
\]
where $q=q_{0}+q_{1}i+q_{2}j+q_{3}k$, $u=u_{0}+u_{1}i+u_{2}j+u_{3}k \in \Sp(1)$, \(c_{\ell}=a_{\ell}+b_{\ell}i\in\mathbb C^{+}\), \(\ell=1,2\), and $0\le \theta\le \frac{\pi}{2}$. 
Then
\begin{align}\label{xAx}
	x^{*}Ax=c_{1}(1-t)+tu^{*}c_{2} u + \rho qu,
\end{align}
where $t=\sin^{2}\theta$, $\rho=\sqrt{t(1-t)}=\sin\theta\cos\theta$. Therefore,
\begin{align}\label{upperbildset}
	\bild^{+}(A)=\{ \pi ( c_{1}(1-t)+tu^{*}c_{2} u + \rho qu) : 0\leq t\leq 1, u\in \Sp(1)\}. 
\end{align}

\begin{lemma}[\cite{K19}]\label{connected}
	Let \(A \in M_{2}(\H) \), and let \(L\) be a line parallel to the
	\(y\)-axis. Then \(\bild^{+}(A) \cap L\) is either empty or
	connected.
\end{lemma}

\begin{remark}
	Lemma \ref{connected} implies that, once the upper and lower boundary points are
	identified for a fixed real part \(X\), the entire vertical segment
	between them belongs to \(\bild^+(A)\). Hence, in the
	following sections, after determining the upper and lower boundary
	curves, we will not separately verify that the region between them is
	filled.
\end{remark}

\begin{remark}
	The shapes of the upper bilds in this section were first determined by So and Thompson \cite[Lemmas~6.1--6.7, Theorem~7.1]{ST1}. Their description is given after the Hermitian part has been diagonalized and normalized to $\operatorname{diag}(0,1)$, so that their parameters are the diagonal entries of the skew-Hermitian part in this basis rather than spectral data of~$A$. In contrast, we work directly with the Schur form of Lemma~\ref{uptri}, whose diagonal entries are the standard eigenvalues $c_1,c_2$ of~$A$. In these parameters the boundary ellipses have the standard eigenvalues (or $c_1$ and $\overline{c_2}$) as foci, the threshold between the cases in Section~\ref{m>0} is $m^2=4b_1b_2$, and the existence of the common tangent segment in Section~\ref{gamma>0} is decided by the sign of~$\Lambda$. 
\end{remark}

\subsection{$q=0$, $c_{1}$ and $c_{2}$ are complex numbers}

\begin{lemma}[\cite{AY95,qnormal,K19}]\label{normal}
	Let \(A\in M_2(\H)\) be normal, with standard eigenvalues
	\[
	c_1=a_1+b_1i,
	\qquad
	c_2=a_2+b_2i,
	\qquad
	b_1,b_2\geq0.
	\]
	If \(b_1+b_2>0\), define
	\[
	v=\frac{b_2a_1+b_1a_2}{b_1+b_2}.
	\]
	Then
	\[
	\bild^+(A)=\operatorname{conv}\{c_1,c_2,v\}.
	\]
	If \(b_1=b_2=0\), then
	\[
	\bild^+(A)=W_{\H}(A)
	=
	[\min\{a_1,a_2\},\max\{a_1,a_2\}].
	\]
\end{lemma}



	
	\begin{proof}
		Since $A$ is normal, there exists a quaternionic unitary matrix $U$ such that
		\[
		U^*AU
		=
		\begin{pmatrix}
			c_1 & 0\\
			0 & c_2
		\end{pmatrix},
		\]
		where
		\[
		c_1=a_1+b_1 i,
		\qquad
		c_2=a_2+b_2 i,
		\qquad
		b_1,b_2\ge 0.
		\]
        
		By (\ref{xAx}), we can see that
		\[
		x^*Ax
		=
		(1-t)c_1+t\,u^*c_2u = (1-t)(a_{1}+b_{1}i)+t(a_2+b_2 u^{*}iu).
		\]
		Then the upper
		complex representative of \(x^*Ax\) is
		\[
		\pi(x^*Ax)
		=
		(1-t)a_1+ta_2
		+
		i\left|(1-t)b_1i+tb_2u^{*}iu\right|.
		\]
		Let \(\phi\in[0,\pi]\) denote the angle between \(i\) and \(u^*iu\).
		Since both \(i\) and \(u^*iu\) are unit pure quaternions, we have
		\[
		\langle i,u^*iu\rangle=\cos\phi.
		\]
		Then
		\[
		\pi(x^*Ax)
		=
		(1-t)a_1+ta_2
		+
		i\sqrt{
			(1-t)^2b_1^2+t^2b_2^2
			+2t(1-t)b_1b_2\cos\phi
		}.
		\]
		Since \(\cos\phi\) ranges over \([-1,1]\), for each fixed
		\(t\in[0,1]\), the imaginary part ranges over
		\[
		\left[
		\left|(1-t)b_1-tb_2\right|,
		(1-t)b_1+tb_2
		\right].
		\]
		Thus, for each fixed \(t\), the corresponding points fill the
		vertical line segment between the upper boundary
		\[
		(1-t)c_1+tc_2
		\]
		and the lower boundary
		\[
		(1-t)a_1+ta_2
		+i\left|(1-t)b_1-tb_2\right|.
		\]
		As \(t\) varies from \(0\) to \(1\), the upper endpoints trace
		\([c_1,c_2]\), while the lower endpoints trace
		\([c_1,v]\cup[v,c_2]\).
	    Therefore,
		\[
		\bild^+(A)=\operatorname{conv}\{c_1,c_2,v\}.
		\]
	\end{proof}
	
	\begin{remark}\label{re=}
		If \(a_1=a_2=a\) and \(b_1+b_2>0\), then
		\(
		v=a,
		\) and 
		\[
		\bild^+(A)
		=
		\{a+it:0\leq t\leq \max\{b_1,b_2\}\}.
		\]
		The full bild is therefore the vertical line segment
		\[
		W_{\H}(A)\cap\C
		=
		\{a+it:|t|\leq\max\{b_1,b_2\}\}.
		\]
		By Theorem~3.1 of~\cite{CDM20star}, the quaternionic numerical range
		is convex if and only if its bild is convex. Hence \(W_{\H}(A)\) is convex.
	\end{remark}
	
	\begin{remark}
		If \(b_1=b_2=0\), then Lemma~\ref{normal} gives
		\[
		W_{\H}(A)
		=
		[\min\{a_1,a_2\},\max\{a_1,a_2\}],
		\]
		which is convex. Thus the two convex cases obtained here,
		\(a_1=a_2\) or \(b_1=b_2=0\), are the two cases in
		Theorem~3.1 of~\cite{So2023}.
		Conversely, if $a_1\neq a_2$ and $b_1+b_2>0$, choose $\ell$ with $b_\ell>0$. Then $a_\ell\neq v$, and the vertical section of $\bild^{+}(A)=\operatorname{conv}\{c_1,c_2,v\}$ at $X=a_\ell$ is $\{c_\ell\}$. Hence the vertical section of the bild at $X=a_\ell$ is $\{a_\ell\pm ib_\ell\}$, which is not convex, and $W_{\H}(A)$ is not convex by Theorem~3.1 of~\cite{CDM20star}. Therefore, the two convex cases above are the only ones, in agreement with Theorem~3.1 of~\cite{So2023}.
	\end{remark}
	
	\subsection{$c_{1},c_{2}$ are real numbers, $q$ is a nonzero quaternion}
	~
	
    In this case, $a_{1}, a_{2}\in \mathbb{R}$, $b_{1}=b_{2}=0$, and $q\in \H \setminus \{0\}$. There is a quaternionic unitary matrix, such that the original matrix is unitarily similar to 
		\[
	A=
	\begin{pmatrix}
		a_1 & |q|\\
		0 & a_2
	\end{pmatrix}.
	\]
    By (\ref{xAx}),
    \[
    x^*Ax=(1-t)a_1+ta_2+|q|\sqrt{t(1-t)}\,u.
    \]
    Hence, writing \(z=X+iY\in\bild^+(A)\),
    \[
    X=(1-t)a_1+ta_2+|q|\sqrt{t(1-t)}\,u_0,
    \qquad
    Y=|q|\sqrt{t(1-t)}\sqrt{1-u_0^2}.
    \]
    Therefore, for fixed \(t\),
    \[
    \left(X-((1-t)a_1+ta_2)\right)^2+Y^2
    =
    |q|^2t(1-t),
    \qquad Y\ge0.
    \]
	Thus
	\[
	\bild^+(A)
	=
	\bigcup_{0\le t\le1}
	\left\{
	X+iY:
	\left(X-\big((1-t)a_1+ta_2\big)\right)^2+Y^2
	=
	|q|^2t(1-t),
	\quad
	Y\ge0
	\right\}.
	\]

	Although each fixed \(t\) gives only a semicircle, the union over
	\(0\le t\le1\) fills an elliptical region, as shown below.
	
	Set
	\[
	c=\frac{a_1+a_2}{2},
	\qquad
	d=a_2-a_1,
	\qquad
	r=t-\frac12,
	\qquad
	P=X-c.
	\]
	Then
	\[
	(1-t)a_1+ta_2=c+dr,
	\qquad
	t(1-t)=\frac14-r^2,
	\]
	so the semicircle equation becomes
	\[
	(P-dr)^2+Y^2
	=
	|q|^2\left(\frac14-r^2\right),
	\]
	or equivalently
	\[
	F(r):=(d^2+|q|^2)r^2-2dPr+P^2+Y^2-\frac{|q|^2}{4}=0.
	\]
	
	For a fixed point \((X,Y)\) with \( Y\geq 0\), membership in the union of semicircles is
	equivalent to the existence of \(r\in[-1/2,1/2]\) such that \(F(r)=0\),
	or equivalently,
	\[
	\min_{-1/2\le r\le1/2}F(r)\le0.
	\]
	At the endpoints,
	\[
	F\!\left(-\frac12\right)
	=
	\left(P+\frac d2\right)^2+Y^2,
	\qquad
	F\!\left(\frac12\right)
	=
	\left(P-\frac d2\right)^2+Y^2.
	\]
	Hence
	\[
	F\!\left(-\frac12\right)=0
	\iff
	(X,Y)=(a_1,0),
	\qquad
	F\!\left(\frac12\right)=0
	\iff
	(X,Y)=(a_2,0).
	\]
    Except at the endpoint points \((a_1,0)\) and \((a_2,0)\),
    \(
    F\!\left(-\frac12\right)>0,
    F\!\left(\frac12\right)>0.
    \)
    Hence, for every other point in the union, the strict convexity of
    \(F\) implies that its unique minimizer
    \[
    r_0=\frac{dP}{d^2+|q|^2}
    \]
    lies in \((-1/2,1/2)\). Substituting \(r_0\) gives
	\[
	F(r_0)
	=
	\frac{|q|^2}{d^2+|q|^2}P^2+Y^2-\frac{|q|^2}{4},
	\]
	so \(F(r_0)\le0\) is equivalent to
	\[
	\frac{P^2}{(d^2+|q|^2)/4}
	+
	\frac{Y^2}{|q|^2/4}
	\le1.
	\]
	Conversely, if $(X,Y)$ satisfies this inequality, then $|P|\le\frac12\sqrt{d^2+|q|^2}$, so
	\(
	|r_0|\le\frac{|d|}{2\sqrt{d^2+|q|^2}}<\frac12 .
	\)
	Since $F(\pm\frac12)\ge0$ and $F(r_0)\le0$, the intermediate value theorem gives some $r\in[-\frac12,\frac12]$ with $F(r)=0$, that is, $(X,Y)\in\bild^{+}(A)$.
	Substituting these quantities back, we get
	\[
	\frac{\left(X-\frac{a_1+a_2}{2}\right)^2}
	{\frac{(a_2-a_1)^2+|q|^2}{4}}
	+
	\frac{Y^2}{\frac{|q|^2}{4}}
	\le 1,
	\qquad Y\ge 0.
	\]
	Thus \(\bild^+(A)\) is the upper half of an elliptical disk.
	
	Reflecting \(\bild^+(A)\) across the real axis shows that the full bild is
	an elliptical disk, and hence is convex. By Theorem~3.1
	of~\cite{CDM20star}, \(W_{\mathbb H}(A)\) is convex.

	\begin{proposition}\label{tri}
		If \[
		A=
		\begin{pmatrix}
			a_1 & q\\
			0 & a_2
		\end{pmatrix},
		\qquad 
		a_1,a_2\in\mathbb R,\quad q\in\mathbb H\setminus \{ 0 \}.
		\]
		Then, $\bild^{+}(A)$ is a half elliptical disk, centered at 
		\(
		\left(\frac{a_1+a_2}{2},0\right),
		\)
		symmetric about the line
		\(
		X=\frac{a_1+a_2}{2}.
		\)
		Its maximum imaginary part is
		\(
		Y_{\max}=\frac{|q|}{2}
		\)
		at
		\(
		X=\frac{a_1+a_2}{2}
		\). It is convex,
		and
		$W_{\mathbb H}(A)$  is convex.
	\end{proposition}

   \begin{remark}
   	Consider
   	\[
   	A=
   	\begin{pmatrix}
   		0&q\\
   		0&0
   	\end{pmatrix},
   	\qquad q\in\mathbb H\setminus \{0\}. 
   	\]
   	Here \(a_1=a_2=0\), so the elliptic formula reduces to
   	\[
   	X^2+Y^2\le \frac{|q|^2}{4},
   	\qquad Y\ge0.
   	\]
   	Hence \(\bild^+(A)\) is the upper half disk centered at the origin with radius
   	\(|q|/2\). Reflecting across the real axis gives the full bild
   	\[
   	\left\{
   	X+iY\in\mathbb C:
   	X^2+Y^2\le \frac{|q|^2}{4}
   	\right\}.
   	\]
   	Thus \(W_{\mathbb H}(A)\) is convex. In particular, our result recovers Case (3) in \cite{K19} as a special case.
   \end{remark}

	\subsection{$c_{1}$ is a real number, $c_{2}$ is a complex number, $q$ is quaternion}
	\
	
	The case in which the nonreal standard eigenvalue appears first and the real one second gives the same upper bild, by taking adjoints and permuting the basis vectors; see Remark~5 of the Supplementary Material. Hence we may assume that $c_1$ is real. In this case, $a_1,a_2,b_2\in\mathbb R, ~ q\in\mathbb H\setminus\{ 0 \}, ~b_{2}> 0$. This is the last case in which the off-diagonal entry can be treated
	through its norm alone. It can be verified that the original matrix is unitarily similar to 
	\[
	A=\begin{pmatrix}
	a_{1} & |q| \\ 0 &a_{2}+b_{2}i
	\end{pmatrix}.
	\]
    By (\ref{xAx}), we have
	\[
	x^{*}Ax
	=
	(1-t)a_1+ta_2
	+
	tb_2\, u^{*}iu
	+
	|q|\sqrt{t(1-t)}\,u.
	\]
	Let \(z=X+iY\in \bild^+(A)\).
	Then
	\[
	X=(1-t)a_1+ta_2
	+
	|q|\sqrt{t(1-t)}\,u_0,  \quad Y
	=
	\left|
	tb_2\, u^{*}iu 
	+
	|q|\sqrt{t(1-t)} \im u
	\right|.
	\]
	Using
	\[
	| u^{*}iu|=1,~ \left\langle u^{*}iu,\im u\right\rangle=u_1,
	\]
	we obtain
	\[
	Y^2
	=
	t^2b_2^2
	+
	|q|^2t(1-t)(1-u_0^2)
	+
	2tb_2 |q|\sqrt{t(1-t)}\,u_1.
	\]
	For fixed \(t\) and feasible \(X\), the formula for \(X\) fixes
	\(u_0\). Since \(u\) is a unit quaternion,
	\[
	u_1^2+u_2^2+u_3^2=1-u_0^2.
	\]
	Thus \(u_1\) ranges over
	\([-\sqrt{1-u_0^2},\sqrt{1-u_0^2}]\). Every value in this interval
	is attainable by a suitable choice of \(u_2,u_3\).
	The displayed formula for \(Y^2\) is increasing in \(u_1\), so the
	maximum and minimum attainable heights for these fixed \(t,X\) are
	\[
	Y_{\max}
	=tb_2+|q|\sqrt{t(1-t)}\sqrt{1-u_0^2},
	\qquad
	Y_{\min}
	=\left|tb_2-|q|\sqrt{t(1-t)}\sqrt{1-u_0^2}\right|.
	\]

	Together with
	\begin{equation}\label{rcqx}
		X=(1-t)a_1+ta_2
		+
		|q|\sqrt{t(1-t)}\,u_0,
	\end{equation}
	this gives
	\begin{align}\label{rcc}
	\left(
	X-\big((1-t)a_1+ta_2\big)
	\right)^2
	+
	\left(
	Y_{\max}-tb_2
	\right)^2
	=
	|q|^2t(1-t).
   \end{align}
	Thus, for each fixed \(t\), the upper boundary of the corresponding
	slice lies on the upper semicircle centered at
	$
	(1-t)a_1+t(a_2+b_2i)
	$
	with radius
	\(
	|q|\sqrt{t(1-t)}.
	\)
	
    To identify the geometry of the upper envelope, introduce the
	auxiliary filled disks
	\[
	D_t=\left\{(X,Y)\in\mathbb R^2:
	\left(X-\bigl((1-t)a_1+ta_2\bigr)\right)^2
	+(Y-tb_2)^2\le |q|^2t(1-t)\right\},
	\qquad 0\le t\le1,
	\]
	and let \(E=\bigcup_{0\le t\le1}D_t\).
	For fixed \(t\) and feasible \(X\), the highest point of \(D_t\)
	is precisely the attainable value \(Y_{\max}\) above. The interiors
	of the disks are auxiliary and are not asserted to be attainable
	for each fixed \(t\).
	
	Set
	\[
	C=\left(\frac{a_1+a_2}{2},\frac{b_2}{2}\right),
	\quad
	\Delta=a_2-a_1,
	\quad
	v=(\Delta,b_2),
	\quad
	|v|=\sqrt{\Delta^2+b_2^2},
    \quad
	r=t-\frac12,
	\quad
	-\frac12\le r\le\frac12.
	\]
	Then
	\[
	\bigl((1-t)a_1+ta_2,\;tb_2\bigr)=C+rv,
	\qquad
	t(1-t)=\frac14-r^2.
	\]
	For an arbitrary point \((X,Y)\), set \(P=(X,Y)-C\).
	By the definition of \(E\), we have \((X,Y)\in E\) if and only if
	there exists \(r\in[-\frac12,\frac12]\) such that
	\[
	|P-rv|^2\le |q|^2\left(\frac14-r^2\right).
	\]
	Equivalently, there exists such an \(r\) with
	\[
	F(r)\le0,
	\qquad
	F(r):=\bigl(|v|^2+|q|^2\bigr)r^2
	-2(P\cdot v)r+|P|^2-\frac{|q|^2}{4}.
	\]
	Thus \((X,Y)\in E\) is equivalent to
	\[
	\min_{-\frac12\le r\le\frac12}F(r)\le0.
	\]

	At the endpoints,
	\[
	F\!\left(-\frac12\right)
	=
	\left|P+\frac12v\right|^2,
	\qquad
	F\!\left(\frac12\right)
	=
	\left|P-\frac12v\right|^2.
	\]
	Thus both endpoint values are nonnegative, with equality only at the
	two degenerate endpoint points
	\(
	(a_1,0)
    \) and
    \(
	(a_2,b_2).
	\)
	For every other point under
	consideration,
	\[
	F\!\left(-\frac12\right)>0,
	\qquad
	F\!\left(\frac12\right)>0.
	\]
	Since \(F\) is strictly convex, whenever its minimum on
	\([-\frac12,\frac12]\) is nonpositive, the unique unconstrained
	minimizer
	\[
	r_0=\frac{P\cdot v}{|v|^2+|q|^2}
	\]
	must lie in
	\[
	\left(-\frac12,\frac12\right).
	\]
	Hence the constrained minimum is attained at \(r_0\), and
	\[
	F(r_0)
	=
	|P|^2-\frac{|q|^2}{4}
	-\frac{(P\cdot v)^2}{|v|^2+|q|^2}
	\le0.
	\]
	Decompose \(P\) into the directions parallel and perpendicular to \(v\):
	\[
	\xi=\frac{P\cdot v}{|v|},
	\qquad
	\zeta=\frac{P\cdot v^\perp}{|v|},
	\qquad
	v^\perp=(-b_2,\Delta).
	\]
	Since
	\[
	|P|^2=\xi^2+\zeta^2,
	\qquad
	P\cdot v=|v|\xi,
	\]
	the inequality becomes
	\[
	\frac{|q|^2}{|v|^2+|q|^2}\xi^2+\zeta^2
	\le \frac{|q|^2}{4},
	\]
	or equivalently,
	\[
	\frac{\xi^2}{(|v|^2+|q|^2)/4}
	+
	\frac{\zeta^2}{|q|^2/4}
	\le1.
	\]
	Conversely, every point satisfying this elliptical inequality has
	\[
	|\xi|\le\frac12\sqrt{|v|^2+|q|^2},
	\qquad
	|r_0|
	=\frac{|v|\,|\xi|}{|v|^2+|q|^2}
	\le\frac{|v|}{2\sqrt{|v|^2+|q|^2}}
	<\frac12.
	\]
	Thus \(r_0\) is admissible, and the preceding inequality describes
	precisely the auxiliary region \(E\). It is an elliptical disk with
	center \(C\), semiaxes
	\[
	\frac12\sqrt{|v|^2+|q|^2}
	\quad\text{and}\quad
	\frac{|q|}{2},
	\]
	and foci
	\[
	C-\frac v2=(a_1,0),
	\qquad
	C+\frac v2=(a_2,b_2).
	\]
	For each feasible \(X\), the highest point of \(E\) is obtained by
	maximizing the highest points of the disks \(D_t\) over admissible
	\(t\). Each of these highest points is the attainable value
	\(Y_{\max}\) computed above. Hence the upper boundary of
	\(\bild^+(A)\) is the upper elliptic arc of \(E\).

  We now identify the lower boundary. For a fixed feasible \(X\), let
  \[
  I_X
  =
  \left\{
  t\in[0,1]:
  |q|^2t(1-t)
  -
  \left(X-\bigl((1-t)a_1+ta_2\bigr)\right)^2
  \geq0
  \right\}.
  \]
The expression defining \(I_X\) is a strictly concave quadratic
function of \(t\). Since \(X\) is feasible, \(I_X\) is a nonempty
compact interval, possibly a singleton. For \(t\in I_X\), put
  \[
  \phi_X(t)
  =
  tb_2-
  \sqrt{
  	|q|^2t(1-t)
  	-
  	\left(X-\bigl((1-t)a_1+ta_2\bigr)\right)^2
  }.
  \]
At \(t=0\) and \(t=1\), the original parametrization gives
\((X,Y)=(a_1,0)\) and \((a_2,b_2)\), respectively. These endpoint
values are also included in the definitions of \(I_X\) and
\(\phi_X\) above. For each \(t\in I_X\), the lowest point of the auxiliary disk
  \(D_t\) at \(X\) has height \(\phi_X(t)\). Thus the lower endpoint
  of the vertical section of \(E\) is
  \[
  L_E(X)=\min_{t\in I_X}\phi_X(t).
  \]
  By the computation of \(Y_{\min}\) above, the lowest height actually
  attainable in \(\bild^+(A)\) is
  \[
  L_{\bild^+}(X)=\min_{t\in I_X}|\phi_X(t)|.
  \]
  
  Let \(\alpha_X,\beta_X\) be the endpoints of \(I_X\). The square
  root in \(\phi_X\) vanishes at these endpoints, so
  \[
  \phi_X(\alpha_X)=\alpha_Xb_2\geq0,
  \qquad
  \phi_X(\beta_X)=\beta_Xb_2\geq0.
  \]
  If \(L_E(X)\geq0\), then \(\phi_X(t)\geq0\) throughout \(I_X\).
  Hence \(L_{\bild^+}(X)=L_E(X)\). If \(L_E(X)<0\), continuity and
  the nonnegative endpoint values give some \(t_0\in I_X\) with
  \(\phi_X(t_0)=0\). In this case \(L_{\bild^+}(X)=0\). Consequently,
  \[
  L_{\bild^+}(X)=\max\{0,L_E(X)\}.
  \]
  Here \(L_E(X)\) is the height of the lower elliptic arc of \(E\);
  the maximum with zero replaces any portion below the real axis
  by the chord \(E\cap\mathbb R\).
  
  We have therefore identified both endpoints of every nonempty
  vertical section of \(\bild^+(A)\): they are the endpoints of the
  corresponding section of \(E\cap\mathbb C^+\). By
  Lemma~\ref{connected}, each such section of \(\bild^+(A)\) is
  connected. Hence
  \[
  \bild^+(A)=E\cap\mathbb C^+.
  \]
  Geometrically, its upper boundary is the upper elliptic arc; its
  lower boundary consists of the part of the lower elliptic arc above
  the real axis and the chord \(E\cap\mathbb R\).

	\begin{proposition}\label{prop:one-real}
		Let
		\[
		A=
		\begin{pmatrix}
			a_1&q\\
			0&a_2+b_2i
		\end{pmatrix},
		\qquad
		a_1,a_2,b_2\in\mathbb R,\quad b_2>0,\quad q\neq0.
		\]
		Then
		\[
		\bild^+(A)=E\cap\mathbb C^+
		\]
		is convex, where \(E\) is the elliptical disk determined above, centered at
		\[
		C=\left(\frac{a_1+a_2}{2},\frac{b_2}{2}\right).
		\]
		Its lower boundary consists of the part of the lower elliptic arc in
		\(\mathbb C^+\) and the chord \(E\cap\mathbb R\).
	\end{proposition}

	\begin{remark}
		For the matrices in this proposition, \(W_{\H}(A)\) is not convex.
		The leftmost and rightmost points of the elliptical disk \(E\)
		have respective heights
		\[
		\frac{b_2}{2}
		\left(1\mp\frac{\Delta}{\sqrt{\Delta^2+|q|^2}}\right).
		\]
		Both heights are positive, so both points belong to \(\bild^+(A)\).
		Hence, at either extreme real part, the vertical section of
		\(\bild^+(A)\) is a singleton at positive height.
		The corresponding vertical section of the full bild therefore
		consists of two distinct points \(\{X+iY,X-iY\}\), with \(Y>0\),
		and is not convex.
		By \cite[Theorem~3.1]{CDM20star}, \(W_{\H}(A)\) is not convex.
	\end{remark}
	
	\subsection{$c_{1}, c_{2}$ are complex numbers, and $q$ is a nonzero quaternion}
	
	\
	
	When \(b_1,b_2>0\), the off-diagonal entry \(q\) cannot in general be reduced to \(|q|\) while keeping both diagonal entries fixed. 
	However, there is a diagonal unitary matrix $U$ such that
	
	\[
	   U^{*} \begin{pmatrix}
		   a_1+b_1i&q_{0}+q_{1}i+q_{2}j+q_{3}k\\
		   0&a_2+b_2i
 	   \end{pmatrix} U = \begin{pmatrix}
 	   a_1+b_1 i & m+\gamma j\\
 	   0 & a_2+b_2 i \end{pmatrix}.
	\]
	(see Lemma~3 of the Supplementary Material). Therefore, it is enough to consider
	\[
	A=\begin{pmatrix}
		 a_1+b_1 i & m+\gamma j\\
		0 & a_2+b_2 i
	\end{pmatrix}, \text{ where } m = \sqrt{q_{0}^{2} + q_{1}^{2}},~ \gamma = \sqrt{q_{2}^{2} + q_{3}^{2}}.
	\]
	
	In this section, we mainly treat two special cases, $m>0, \gamma = 0$ and $m=0, \gamma>0$. When $m>0$ and $\gamma>0$, the matrix falls into the generic case of Section~\ref{Section:Generic}; see Lemma~4 of the Supplementary Material for more detail.

	\subsubsection{$m>0, \gamma=0$}\label{m>0}
	
	\
	
    In this subsection, we derive an explicit criterion that determines the
   shape of the upper bild of
   \[
   A=
  \begin{pmatrix}
	a_1+b_1i & m\\
	0 & a_2+b_2i
  \end{pmatrix},
  \qquad m>0.
   \]
   When $m^{2}>4b_{1}b_{2}$, the boundary is an ellipse cut by the real axis. When $m^{2}=4b_{1}b_{2}$, the boundary is an ellipse. When $m^{2}<4b_{1}b_{2}$, it is composed of part of an ellipse and two straight-lines segments. 
   
   By (\ref{xAx}), 
   \[
   x^*Ax
   =
   (1-t)a_1+ta_2
   +
   (1-t)b_1i
   +
   tb_2u^*iu
   +
   m\sqrt{t(1-t)}\,u.
   \]
   Thus, writing \(X+iY\in\bild^+(A)\), we have
	\begin{equation}\label{mx}
		 X=(1-t)a_{1}+ta_{2}+m\sqrt{t(1-t)}u_{0},
	\end{equation}
	\[
	Y=|(1-t)b_{1} i+ tb_{2}u^{*}iu+m\sqrt{t(1-t)} \im u|.
	\]
	Set
	\[
	A_t=(1-t)b_1,\qquad B_t=tb_2,\qquad \rho=\sqrt{t(1-t)}.
	\]
	Using
	\[
	|\operatorname{Im}u|^2=1-u_0^2, ~
	\langle i, \im u\rangle=u_1, ~
	\langle u^{*}iu, \im u\rangle=u_1, ~
	\text{ and }
	\langle i, u^{*}iu\rangle
	=
	1-2(u_2^2+u_3^2),
	\]
	we have that 
	\begin{align*}
			Y^2=
		           A_t^2+B_t^2
		           +
		           m^2\rho^2(1-u_0^2)
		           +
		           2A_tB_t\bigl(1-2(u_2^2+u_3^2)\bigr)
		           +
		           2m\rho(A_t+B_t)u_1.
	\end{align*}

	When \(t=0\) or \(t=1\), we have \(\rho=0\), yielding the endpoint
	values associated with \(a_1+b_1i\) and \(a_2+b_2i\), respectively.
	Thus, for the remaining analysis, we assume \(0<t<1\). For fixed \(t\) and a feasible \(X\), equation~(\ref{mx}) gives
	\[
	u_0=\frac{X-((1-t)a_1+ta_2)}{m\rho},
	\qquad |u_0|\leq1,
	\qquad
	S=\sqrt{1-u_0^2}.
	\]
	Then
	\[
	u_1^2+u_2^2+u_3^2=S^2,
	\qquad
	u_2^2+u_3^2=S^2-u_1^2,
	\]
	and hence
	\begin{equation}\label{phitxu1}
	Y^2:=\Phi_{t,X}(u_1)=
	(A_t+B_t)^2
	+(m^2\rho^2-4A_tB_t)S^2
	+4A_tB_tu_1^2
	+2m\rho(A_t+B_t)u_1,
	\qquad
	-S\le u_1\le S.
	\end{equation}
   
    For fixed \(t\) and \(X\), the extrema of \(Y\) on the corresponding
    slice are determined by optimizing \(\Phi_{t,X}\) over \([-S,S]\).
    We can see that $\Phi_{t,X}(u_1)$
    is a strictly convex quadratic function of \(u_1\) since \( A_t>0,
    B_t>0 \) when $0<t<1$. Moreover,
    \[
    \Phi_{t,X}(S)-\Phi_{t,X}(-S)
    =
    4m\rho(A_t+B_t)S
    \ge0.
    \]
    Therefore the maximum is attained at \(u_1=S\), whereas the minimum is
    attained either at \(u_1=-S\) or at the interior critical point
    \[
    u_1^*
    =
    -\frac{m\rho(A_t+B_t)}{4A_tB_t} \text{ when } u_1^*\in (-S,S).
    \]

	At the endpoint \(u_1=S\),
	\[
	\Phi_{t,X}(S)
	=
	(A_t+B_t+m\rho S)^2.
	\]
	Since \(A_t+B_t>0\), we have
	\[
	Y_{\max}=A_t+B_t+m\rho S=(1-t)b_1+tb_2+m\rho S.
	\]
	Using
	\[
	X-((1-t)a_1+ta_2)=m\rho u_0,
	\qquad
	Y_{\max}-((1-t)b_1+tb_2)=m\rho S,
	\]
	together with \(u_0^2+S^2=1\), we obtain
	\begin{align}\label{circlem}
	\left(X-((1-t)a_1+ta_2)\right)^2
	+
	\left(Y_{\max}-((1-t)b_1+tb_2)\right)^2
	=
	m^2\rho^2 = m^2t(1-t).
    \end{align} 
    For each fixed \(t\), this is the upper semicircle centered at
    \(
    \bigl((1-t)a_1+ta_2,\,(1-t)b_1+tb_2\bigr).
    \)
    Applying the same auxiliary filled-disk calculation as in the preceding
    subsection, we find that the upper boundary is the upper arc of an
    ellipse \(E\) with foci
    \[
    F_1=(a_1,b_1),\qquad F_2=(a_2,b_2),
    \]
    center
    \[
    C=\left(\frac{a_1+a_2}{2},\frac{b_1+b_2}{2}\right),
    \]
    and semi-axis lengths
    \[
    \frac12\sqrt{(a_2-a_1)^2+(b_2-b_1)^2+m^2}
    \qquad\text{and}\qquad
    \frac m2.
    \]

	Next, take \(u_1=-S\). Then
	\[
	\Phi_{t,X}(-S)
	=
	(A_t+B_t-m\rho S)^2,
	\]
	so that
	\[
	Y_{-}=|A_t+B_t-m\rho S|.
	\]
    As in the preceding subsection, minimizing the signed endpoint height
    \(A_t+B_t-m\rho S\) over feasible \(t\), for each fixed \(X\),
    gives the lower boundary of the associated auxiliary elliptical disk \(E\).

We now determine when the critical point \(u_1^*\)
is admissible. Since \(u_1^*<0\), set
\[
g(t):=|u_1^*|
=
\frac{m((1-t)b_1+tb_2)}
{4\sqrt{t(1-t)}\,b_1b_2}.
\]
Interior admissibility requires \(g(t)<S\le1\), and hence \(g(t)<1\). By the arithmetic--geometric mean inequality,
\[
g(t)\geq
\frac{m}{2\sqrt{b_1b_2}},
\]
with the minimum attained at
\[
t_0=\frac{b_1}{b_1+b_2}.
\]

	\subsubsection*{Case 1: \(m^2>4b_1b_2\)}
	
	Then
	\[
	\min_{0<t<1} g(t)
	=
	\frac{m}{2\sqrt{b_1b_2}}
	>1.
	\]
	Hence, for every \(0<t<1\),
	\[
	|u_1^*|=g(t)>1\ge S,
	\]
	so \(u_1^*\notin[-S,S]\). Therefore the interior critical point never contributes to the lower boundary. The lower boundary is determined by the endpoint
	\(u_1=-S\) after upper-bild folding: it consists of the portion of
	the lower elliptic arc above the real axis and the real-axis segment
	cut out by the ellipse. The real-axis segment is attained by the same continuity argument
	as in the proof of Proposition~\ref{prop:one-real}. 
	
	\subsubsection*{Case 2: \(m^2=4b_1b_2\)}
	
	In this case, \(m=2\sqrt{b_1b_2}.\) 
	Therefore
	\[
	\min_{0<t<1}g(t)=g(t_{0})=1, ~t_0=\frac{b_1}{b_1+b_2}.
	\]
	Moreover, the corresponding real-axis point is
	\[
	T=(1-t_0)a_1+t_0a_2
	=
	\frac{b_2a_1+b_1a_2}{b_1+b_2}.
	\]
	At \(X=T\), we have
	\(
	u_0=0
	\)
	,
	\( S=1 \).
	Thus
	\(
	u_1^*=-1=-S.
	\)
	Therefore, in the critical case
	\(
	m^2=4b_1b_2,
	\)
	the interior critical point degenerates to the endpoint \(u_1=-S\). It does
	not give a new boundary piece. Geometrically, this corresponds to the
	associated ellipse being tangent to the real axis. Therefore, together with Lemma $\ref{connected}$, the upper bild is convex. 
	
	\begin{remark}
		The threshold \(m^2=4b_1b_2\) has a direct geometric interpretation.
		Recall that the signed lower branch of the associated ellipse is
		\[
		Y_{\rm ell,lower}
		=
		(1-t)b_1+tb_2-m\sqrt{t(1-t)}\,S.
		\]
		Since \(S\le1\),
		\[
		Y_{\rm ell,lower}
		\ge
		(1-t)b_1+tb_2-m\sqrt{t(1-t)}.
		\]
		When \(m^2=4b_1b_2\), we have
		\[
		Y_{\rm ell,lower}
		\ge
		\left(
		\sqrt{(1-t)b_1}-\sqrt{tb_2}
		\right)^2
		\ge0,
		\]
		with equality only at
		\[
		t=t_0=\frac{b_1}{b_1+b_2},
		\qquad
		S=1.
		\]
		Hence the associated ellipse is tangent to the real axis at
		\[
		(T,0),
		\qquad
		T=\frac{b_2a_1+b_1a_2}{b_1+b_2}.
		\]
		If \(m^2>4b_1b_2\), then at \(t=t_0\) and \(S=1\),
		\(
		Y_{\rm ell,lower}<0.
		\)
		Thus the lower branch crosses the real axis, and consequently the
		intersection of the associated ellipse with the real axis is a
		nontrivial interval.
	\end{remark}
	\subsubsection*{Case 3: \(m^2<4b_1b_2\)}
	
   In this case,
   \[
   g(t_0)
   =
   \frac{m}{2\sqrt{b_1b_2}}
   <1,
   \qquad
   t_0=\frac{b_1}{b_1+b_2}.
   \]
   Set
   \[
   T=(1-t_0)a_1+t_0a_2
   =
   \frac{b_2a_1+b_1a_2}{b_1+b_2}.
   \]
  At \(t=t_0\) and \(X=T\), we have \(u_0=0\) and \(S=1\). Hence
  \[
  |u_1^*(t_0)|=g(t_0)<S.
  \]
  Thus \(u_1^*\) is a genuine interior critical point at \((T,t_0)\).
  By continuity, it remains admissible in a neighborhood of \((T,t_0)\).

   Substituting \(u_1=u_1^*\) into~\eqref{phitxu1} and using
   \[
   A_tB_t=\rho^2b_1b_2,
   \qquad
   S^2=1-u_0^2,
   \]
   we obtain
   \[
   \Phi_{t,X}(u_1^*)
   =
   (4b_1b_2-m^2)
   \left[
   \rho^2u_0^2
   +
   \frac{(A_t-B_t)^2}{4b_1b_2}
   \right].
   \]

    Since
    \[
    A_t-B_t=(1-t)b_1-tb_2,  ~~  \rho u_0
    =
    \frac{
    	X-\bigl((1-t)a_1+ta_2\bigr)
    }{m},
    \]
    it follows that
    \begin{align}\label{phiu_{1}*}
    \Phi_{t,X}(u_1^*)
    =
    (4b_1b_2-m^2)
    \left[
    \frac{
    	\left(
    	X-\bigl((1-t)a_1+ta_2\bigr)
    	\right)^2
    }{m^2}
    +
    \frac{
    	\left(
    	(1-t)b_1-tb_2
    	\right)^2
    }{4b_1b_2}
    \right] :=(4b_1b_2-m^2)F_X(t),
    \end{align}
    where
    \[
    F_X(t)
    =
    \frac{(X-a_1-dt)^2}{m^2}
    +
    \frac{(b_1-ct)^2}{4b_1b_2}, ~d=a_2-a_1,~
    c=b_1+b_2.
    \]
    For each fixed \(X\), \(F_X(t)\) is a strictly convex quadratic
    function of \(t\). Differentiating gives
    \[
    F_X'(t)
    =
    -\frac{2d(X-a_1-dt)}{m^2}
    -\frac{c(b_1-ct)}{2b_1b_2}.
    \]
    Hence the unique critical point satisfies
    \[
    4b_1b_2d(X-a_1-dt)
    +
    m^2c(b_1-ct)
    =
    0.
    \]
    Therefore
    \[
    t_X^*
    =
    \frac{
    	4b_1b_2d(X-a_1)+m^2cb_1
    }{
    	4b_1b_2d^2+m^2c^2
    }.
    \]
    
    Let
    \[
    D
    =
    4b_1b_2d^2+m^2c^2.
    \]
    Completing the square, we obtain
   \begin{align} \label{eq:F-square}
    F_X(t)
    =
    \frac{D}{4m^2b_1b_2}
    (t-t_X^*)^2
    +
    \frac{
    	\bigl[c(X-a_1)-db_1\bigr]^2
    }{D}.
    \end{align}
    Since
    \[
    c(X-a_1)-db_1
    =
    (b_1+b_2)X-(b_2a_1+b_1a_2),
    \]
    we obtain
    \[
    \min_t F_X(t)
    =
    \frac{[(b_1+b_2)X-(b_2a_1+b_1a_2)]^2}{D}.
    \]
    Consequently, the formal interior critical-value envelope
    \(Y_{\mathrm{int}}(X)\) satisfies
    \[
    Y_{\mathrm{int}}(X)^2
    =(4b_1b_2-m^2)\min_{t\in\mathbb R}F_X(t)
    =\frac{4b_1b_2-m^2}{D}
    \bigl[(b_1+b_2)X-(b_2a_1+b_1a_2)\bigr]^2.
    \]

    Since
    \[
    T
    =
    \frac{b_2a_1+b_1a_2}{b_1+b_2},
    \]
    we obtain
    \begin{equation*}
    Y_{\mathrm{int}}(X)^2
    =
    \frac{
    	(4b_1b_2-m^2)(b_1+b_2)^2
    }{D}
    (X-T)^2.
    \end{equation*}
    Thus
     \begin{equation}\label{lineq}
    	Y_{\mathrm{int}}(X)=\kappa|X-T|,
    \end{equation}
    where
    \[
    	\kappa
    	=
    	\frac{
    		(b_1+b_2)\sqrt{4b_1b_2-m^2}
    	}{
    		\sqrt{
    			4b_1b_2(a_2-a_1)^2
    			+
    			m^2(b_1+b_2)^2
    		}
    	}.
    \]
    Hence the interior-minimum envelope consists of two straight-line
    branches,
    \[
    Y_{\mathrm{int}}(X)=
    \begin{cases}
    	\kappa(T-X), & X<T,\\[1mm]
    	\kappa(X-T), & X>T,
    \end{cases}
    \]
    restricted to the region in which \(t_X^*\) and \(u_1^*\) are
    admissible. The two branches meet at
    \((T,0)\).
    Notice that
    \[
    t_T^*
    =
    \frac{b_1}{b_1+b_2}
    =
    t_0.
    \]
    Since the interior critical point is admissible at \((T,t_0)\), as
    verified above, the point \((T,0)\) is actually attained.

    For fixed feasible \(X\) and \(0<t<1\), the interior critical point
    \(u_1^*\) is admissible precisely when \(|u_1^*|<S\).
    Since \(\rho^2=t(1-t)>0\) and
    \[
    H(X,t):=\rho^2\bigl(S^2-(u_1^*)^2\bigr),
    \]
    this is equivalent to \(H(X,t)>0\),
    where
    \[
    H(X,t)=t(1-t)-\frac{(X-a_1-dt)^2}{m^2}
    -\frac{m^2((1-t)b_1+tb_2)^2}{16b_1^2b_2^2}.
    \]
    Define
    \[
    J:=\{X\in\mathbb R:H(X,t_X^*)>0\}.
    \]
    Since \(t_X^*\) is affine in \(X\), \(H(X,t_X^*)\) is a
    strictly concave quadratic function of \(X\). It is positive
    at \(X=T\), so \(J\) is a bounded open interval. Write \(J=(X_-,X_+)\). At either endpoint
    \(X_r\in\{X_-,X_+\}\), we have
    \(0<t_{X_r}^*<1\) and
    \(u_1^*(t_{X_r}^*)=-S(X_r,t_{X_r}^*)\).
    Thus \(J\) has exactly two transition endpoints.
    
    For a fixed feasible \(X\), compactness of the original parameter
    fiber ensures that the lowest height is attained. Although \(u_1^*\) need not be admissible, completing the square in
    \eqref{phitxu1} gives, for every feasible \(0<t<1\),
    \begin{equation}\label{eq:Phi-square}
    	\Phi_{t,X}(u_1)
    	=(4b_1b_2-m^2)F_X(t)
    	+4A_tB_t(u_1-u_1^*)^2.
    \end{equation}
    Combining \eqref{eq:Phi-square} with \eqref{eq:F-square} gives
    \[
    \Phi_{t,X}(u_1)
    \ge (4b_1b_2-m^2)F_X(t)
    \ge \kappa^2(X-T)^2.
    \]
    The same bound holds at \(t=0,1\) by continuity.
    For \(X\in J\), equality is attained at
    \(t=t_X^*\), \(u_1=u_1^*(t_X^*)\).
    
    For \(X\notin J\), a strict interior minimizer would require
    \(u_1=u_1^*(t)\) and \(F_X'(t)=0\), hence \(t=t_X^*\),
    contrary to interior inadmissibility.
    By the endpoint analysis above, the minimum must therefore occur
    on \(u_1=-S\), including the limiting endpoint cases.
    Since \(m<2\sqrt{b_1b_2}\), the arithmetic--geometric mean
    inequality gives \(A_t+B_t-m\rho S>0\).
    Thus the endpoint envelope is the lower elliptic arc of \(E\).
    
    We now examine how the straight-line branch joins the endpoint branch. By (\ref{phiu_{1}*}) and (\ref{eq:Phi-square}), we have that
    \[
    \Phi_{t,X}(u_1)
    =
    \Phi_{t,X}(u_1^*)
    +
    4A_tB_t(u_1-u_1^*)^2.
    \]
    In particular, at the lower endpoint \(u_1=-S\),
    \[
    \Phi_{t,X}(-S)-\Phi_{t,X}(u_1^*)
    =
    4A_tB_t(S+u_1^*)^2.
    \]
    Thus, when
    \(
    -S<u_1^*<S,
    \)
    the constrained minimum with respect to \(u_1\) is attained at the
    interior critical point \(u_1^*\), and its envelope is the
    straight-line branch obtained above. When the interior critical point
    reaches the lower endpoint, that is
    \(
    u_1^*=-S,
    \)
    then the two values coincide, that is
    \(
    \Phi_{t,X}(-S)=\Phi_{t,X}(u_1^*).
    \)
    When \(u_1^*\leq-S\), the constrained minimum is attained at
    \(u_1=-S\), whose envelope is the elliptic branch obtained previously.

    It remains to verify that the straight-line branch and the elliptic
    branch meet tangentially. Define
    \[
    I(t,X)=\Phi_{t,X}(u_1^*(t)), ~E(t,X)=\Phi_{t,X}(-S(t,X)).
    \]
    Then
   \begin{equation}\label{df}
   	  E(t,X)-I(t,X)
   	 =
   	 4A_tB_t(S+u_1^*)^2.
   \end{equation} 
    Let $(X_r,t_r)$ be a transition point. Since
    \begin{equation}\label{00}
    S(X_r,t_r)+u_1^*(t_r)=0,
    \end{equation}
    we have
    \(
    E(t_r,X_r)=I(t_r,X_r).
    \)
  Moreover, \eqref{df} implies that
  \[
  E(t,X_r)\geq I(t,X_r)
  \]
  for every admissible \(t\). The straight-line branch is obtained by minimizing $I(t,X)$
    with respect to $t$. Hence, at
    \(
    t_r=t_{X_r}^*,
    \)
    we have
    \(
    I(t,X_r)\geq I(t_r,X_r)
    \)
    for every admissible $t$. Therefore,
    \[
    E(t,X_r)
    \geq I(t,X_r)
    \geq I(t_r,X_r)
    =
    E(t_r,X_r).
    \]
    Thus, $t_r$ also minimizes $E(t,X_r)$. Consequently, at the
    transition point, the straight-line envelope and the elliptic
    envelope are realized by the same parameter $t_r$. Let
    \[
    L(X):=\min_t I(t,X),
    \qquad
    \mathcal E(X):=\min_t E(t,X).
    \]
    At the transition point \(X=X_r\), both minima are attained at
    \(t=t_r\). Hence
    \[
    L(X_r)=I(t_r,X_r),
    \qquad
    \mathcal E(X_r)=E(t_r,X_r).
    \]
    It follows that 
    \[
    \mathcal E(X_r)=E(t_r,X_r)=I(t_r,X_r)=L(X_r).
    \]
    
    Near \(X_r\), the parameter \(t_X^*\) remains feasible, so
    \(L(X)=\kappa^2(X-T)^2\).
    Moreover, \(E(t,X)\ge I(t,X)\) gives
    \(\mathcal E(X)\ge L(X)\), with equality at \(X=X_r\).
    Both functions are differentiable near \(X_r\):
    \(L\) is quadratic, and \(\mathcal E\) is the square of the
    positive, smooth lower elliptic height.
    Hence \(\mathcal E-L\) has a local minimum at \(X_r\), and therefore
    \[
    \mathcal E'(X_r)=L'(X_r).
    \]

    Hence the squared elliptic envelope and the squared straight-line
    envelope have the same value and the same first derivative at the
    transition point. Since the transition point has positive
    \(Y\)-coordinate, the corresponding branches of \(Y\) also have the
    same value and the same first derivative. Therefore the straight-line
    branch and the elliptic branch meet tangentially.
 
\begin{proposition}
	Let
	\[
	A=
	\begin{pmatrix}
		a_1+b_1i&m\\
		0&a_2+b_2i
	\end{pmatrix},
	\qquad b_1,b_2,m>0,
	\]
	and let \(E\) be the elliptical disk determined above. Then
	\(\bild^+(A)\) is convex. In addition, let
	\(
	T=\frac{b_2a_1+b_1a_2}{b_1+b_2}.
	\)
	\begin{enumerate}
		\item If \(m^2>4b_1b_2\), then
		\(
		\bild^+(A)=E\cap\mathbb C^+.
		\)
		
		\item If \(m^2=4b_1b_2\), then \(E\) is tangent to the real axis at
		\((T,0)\), and
		\(
		\bild^+(A)=E.
		\)
		
		\item If \(m^2<4b_1b_2\), the lower boundary consists of two portions
		of the lower elliptic arc and two tangent line segments meeting at
		\((T,0)\). Equivalently, $\bild^{+}(A)=\operatorname{conv}\bigl(E\cup\{(T,0)\}\bigr)$.
	\end{enumerate}
\end{proposition}

\begin{remark}
	For the matrices in this proposition, \(W_{\H}(A)\) is not convex.
	Indeed, at either extreme real part, the upper bild has a singleton
	vertical section of positive height, so the corresponding section
	of the full bild consists of two distinct points.
\end{remark}

\begin{remark}
		For complex matrices, a different, complex-analytic description of
		the upper bild, as the convex hull of the parts of $W_{\C}(A)$ and
		$W_{\C}(A^*)$ in $\C^+$ (elliptical disks in the $2\times2$ case)
		together with real points, was obtained by Carvalho, Diogo and
		Mendes. We refer the interested reader to \cite{CDM21complex}.
\end{remark}

	\subsubsection{The case \(m=0\) and \(\gamma>0\)}\label{gamma>0}

	By the reduction discussed above, we may write
	\(
	q=\gamma j,
	\)
	with
	\(
	 \gamma>0,
	\)
	and consider
	\[
	A=
	\begin{pmatrix}
		a_1+b_1i&\gamma j\\
		0&a_2+b_2i
	\end{pmatrix},
	\qquad
	a_1,a_2\in\mathbb R,\quad b_1,b_2>0.
	\]
	By (\ref{xAx}), we have
	\[
	x^{*}Ax
	=
	(1-t)a_1+ta_2
	+
	(1-t)b_1i
	+
	tb_2 u^{*}iu
	+
	\gamma \rho\,ju.
	\]
	Let $X+iY\in \bild^{+}(A)$. Then
	\[
	X=(1-t)a_1+ta_2-\gamma \rho u_2,
	\]
	and
	\[
	Y=
	\left|
	A_ti+B_t u^{*}iu+\gamma \rho\left(u_3i+u_0j-u_1k\right)
	\right|,
	\]
	since
	\(
	ju=-u_2+u_3i+u_0j-u_1k.
	\)
	Set
	$
	A_t=(1-t)b_1,
	B_t=tb_2.
	$
	A direct calculation gives
	\[
	\begin{aligned}
		Y^2
		=
		A_t^2+B_t^2
		+
		\gamma^2\rho^2(1-u_2^2)
		+		2A_tB_t\bigl(1-2(u_2^2+u_3^2)\bigr)
		+
		2\gamma \rho(A_t-B_t)u_3.
	\end{aligned}
	\]
	Fix \(t\in(0,1)\) and a feasible \(X\), so that
	\[
	\left|X-\bigl((1-t)a_1+ta_2\bigr)\right|\leq\gamma\rho.
	\]
	Then \(\rho=\sqrt{t(1-t)}>0\), and hence
	\[
	u_2=
	\frac{(1-t)a_1+ta_2-X}{\gamma\rho}.
	\]
	At \(t=0,1\), the original parametrization directly gives
	\((a_1,b_1)\) and \((a_2,b_2)\), respectively.
	Set
	\[
	S=\sqrt{1-u_2^2}.
	\]
	Then
	\[
	u_0^2+u_1^2+u_3^2=S^2.
	\]
	Thus, for fixed \(t\) and feasible \(X\), the remaining optimization
	reduces to
	\[
	Y^2:=\Psi_{t,X}(u_3)=	(A_t-B_t)^2+\gamma^2\rho^2S^2
	+4A_tB_t(S^2-u_3^2)
	+2\gamma\rho(A_t-B_t)u_3,
	\qquad
	-S\le u_3\le S.
	\]
	Since \(b_1,b_2>0\) and \(0<t<1\), we have \(A_tB_t>0\). Hence
	\[
	\Psi_{t,X}''(u_3)=-8A_tB_t<0,
	\]
	so \(\Psi_{t,X}\) is strictly concave on \([-S,S]\).
	Therefore its minimum is attained at one of the endpoints
	\(u_3=\pm S\), whereas its maximum is attained at the critical point
	\[
	u_3^*=\frac{\gamma\rho(A_t-B_t)}{4A_tB_t}
	\]
	when \(|u_3^*|<S\), and otherwise at an endpoint. Moreover,
	\[
	\Psi_{t,X}(S)-\Psi_{t,X}(-S)
	=
	4\gamma\rho(A_t-B_t)S.
	\]
     Hence, if \(S>0\),
	\[
	A_t-B_t>0
	\quad\Longrightarrow\quad
	\Psi_{t,X}(S)>\Psi_{t,X}(-S),
	\]
	whereas
	\[
	A_t-B_t<0
	\quad\Longrightarrow\quad
	\Psi_{t,X}(-S)>\Psi_{t,X}(S).
	\]
	Since
	\(
	\operatorname{sgn}(u_3^*)
	=
	\operatorname{sgn}(A_t-B_t),
	\)
	the preceding observations completely determine the extrema. If
	\(A_t-B_t>0\), the minimum is attained at \(u_3=-S\), while the
	maximum is attained at \(u_3^*\) when \(0<u_3^*<S\), and at
	\(u_3=S\) when \(u_3^*\ge S\). Similarly, if \(A_t-B_t<0\), the
	minimum is attained at \(u_3=S\), while the maximum is attained at
	\(u_3^*\) when \(-S<u_3^*<0\), and at \(u_3=-S\) when
	\(u_3^*\le -S\). If \(A_t=B_t\), then \(u_3^*=0\) and
	\(
	\Psi_{t,X}(S)=\Psi_{t,X}(-S).
	\)	For \(S>0\), the two endpoints give the same minimum value, while
	the maximum is attained at \(u_3=0\).

	At the endpoints,
	\[
	\Psi_{t,X}(S)=(A_t-B_t+\gamma\rho S)^2,
	\qquad
	\Psi_{t,X}(-S)=(A_t-B_t-\gamma\rho S)^2.
	\]
	Hence
	\[
	Y_S=|A_t-B_t+\gamma\rho S|,
	\qquad
	Y_{-S}=|A_t-B_t-\gamma\rho S|.
	\]
	
	To distinguish the two endpoint branches after the upper-bild folding,
	we consider the sign of \(A_t-B_t\).
	
	If \(A_t-B_t>0\), then
	\[
	Y_S=A_t-B_t+\gamma\rho S.
	\]
	Together with
	\[
	X=(1-t)a_1+ta_2-\gamma\rho u_2,
	\qquad
	u_2^2+S^2=1,
	\]
	this gives
	\[
	\left(X-\bigl((1-t)a_1+ta_2\bigr)\right)^2
	+
	\left(Y_S-\bigl((1-t)b_1-tb_2\bigr)\right)^2
	=
	\gamma^2t(1-t).
	\]

   Set
   \[
   C_+
   =
   \left(
   \frac{a_1+a_2}{2},
   \frac{b_1-b_2}{2}
   \right),
   \qquad
   v_+
   =
   \bigl(a_2-a_1,-(b_1+b_2)\bigr), \qquad r=t-\frac12.
   \]
   Then
   \[
   t(1-t)=\frac14-r^2,
   \]
   and the auxiliary full-circle family is given by
   \[
   \left|(X,y)-C_+-rv_+\right|^2
   =
   \gamma^2\left(\frac14-r^2\right),
   \qquad
   -\frac12 \leq r \leq \frac12.
   \]
   The formula extends continuously to \(r=\pm\frac12\), corresponding
   to \(t=0,1\), where the circle radius vanishes and the degenerate
   circles reduce to \((a_1,b_1)\) and \((a_2,-b_2)\), respectively.

   Let \(D_t^+\) be the filled disk with center
   \(C_++(t-\tfrac12)v_+\) and radius
   \(\gamma\sqrt{t(1-t)}\), for \(0\leq t\leq1\).
   As in the preceding subsection, the auxiliary union
   \(E_+=\bigcup_{0\leq t\leq1}D_t^+\) is the elliptical disk,
   \[
   \frac{\xi^2}{\bigl(|v_+|^2+\gamma^2\bigr)/4}
   +
   \frac{\zeta^2}{\gamma^2/4}
   \le 1,
   \]
   where \(\xi\) and \(\zeta\) are the coordinates parallel and perpendicular
   to \(v_+\), respectively.
   Hence the semi-major and semi-minor axes are
   \[
   A_{\rm maj}^+
   =
   \frac12\sqrt{|v_+|^2+\gamma^2},
   \qquad
   A_{\rm min}^+
   =
   \frac{\gamma}{2}.
   \]
   Since the major axis is parallel to \(v_+\), the foci are
   \[
   C_+-\frac12v_+=(a_1,b_1),
   \qquad
   C_++\frac12v_+=(a_2,-b_2).
   \]
   
  If \(A_t-B_t<0\), then
  \[
  Y_{-S}=B_t-A_t+\gamma\rho S.
  \]
  Together with
  \[
  X=(1-t)a_1+ta_2-\gamma\rho u_2,
  \qquad
  u_2^2+S^2=1,
  \]
  we obtain
  \[
  \left(X-\bigl((1-t)a_1+ta_2\bigr)\right)^2
  +
  \left(Y_{-S}-\bigl(-(1-t)b_1+tb_2\bigr)\right)^2
  =
  \gamma^2t(1-t).
  \]
  
  The corresponding auxiliary full-circle family  is the reflection across the real axis of the preceding
  one. Hence the associated auxiliary ellipse \(E_-\) is the reflection
  of \(E_+\), with center, major-axis direction, and foci, 
  \[
  C_-=
  \left(
  \frac{a_1+a_2}{2},
  \frac{-b_1+b_2}{2}
  \right),
  \quad
  v_-=(a_2-a_1,b_1+b_2),
  \quad
  (a_1,-b_1),
  \quad
  (a_2,b_2),
  \]
  respectively. 
  Its semi-major and semi-minor axes are the same as those of \(E_+\).

	If $A_{t}-B_{t}<0$, for \(u_3=S\),
	\[
	Y_S=|A_t-B_t+\gamma\rho S|,
	\]
	so that
	\[
	Y_S
	=
	\begin{cases}
		B_t-A_t-\gamma\rho S,
		& \gamma\rho S\le B_t-A_t,\\[1mm]
		A_t-B_t+\gamma\rho S,
		& \gamma\rho S>B_t-A_t.
	\end{cases}
	\]
	If $A_{t}-B_{t}>0$, for \(u_3=-S\),
	\[
	Y_{-S}=|A_t-B_t-\gamma\rho S|,
	\]
	so that
	\[
	Y_{-S}
	=
	\begin{cases}
		A_t-B_t-\gamma\rho S,
		& \gamma\rho S\le A_t-B_t,\\[1mm]
		B_t-A_t+\gamma\rho S,
		& \gamma\rho S>A_t-B_t.
	\end{cases}
	\]
	Thus the folded endpoint branches are contained in the auxiliary
	elliptical disks \(E_+\cup E_-\). Finally, if \(A_t=B_t\), then \(u_3^*=0\) and
	\[
	Y_S=Y_{-S}=\gamma\rho S.
	\]
	
	If \(S=0\), then
	\[
	Y_S=Y_{-S}=|A_t-B_t|.
	\]
     When \(A_t=B_t\) and \(S>0\), the two endpoints give the same
    minimum value, while the maximum is attained at \(u_3=0\).
	Moreover,
	\[
	A_t=B_t
	\quad\Longleftrightarrow\quad
	t=\frac{b_1}{b_1+b_2},
	\]
	which is precisely the transition value at which the larger endpoint
	changes from \(u_3=S\) to \(u_3=-S\).
	Thus, after the upper-bild folding, the endpoint branches
	\(u_3=\pm S\) introduce no geometric components other than the two
	auxiliary ellipses \(E_+\) and \(E_-\), which are reflections of each
	other across the real axis.
	
	We now consider the interior critical branch. Suppose
	\[
	|u_3^*|<S,
	\qquad
	u_3^*=\frac{\gamma\rho(A_t-B_t)}{4A_tB_t}.
	\]
	Then \(u_3^*\) is the maximizer of \(\Psi_{t,X}\) on \([-S,S]\). Substitution gives
	\[
	\Psi_{t,X}(u_3^*)
	=
	(4A_tB_t+\gamma^2\rho^2)
	\left(
	S^2+\frac{(A_t-B_t)^2}{4A_tB_t}
	\right).
	\]
	Since
	\[
	S^2
	=
	1-
	\frac{\bigl(X-(1-t)a_1-ta_2\bigr)^2}
	{\gamma^2\rho^2},
	\]
	we have
	\[
	\rho^2S^2
	=
	\rho^2-
	\frac{\bigl(X-(1-t)a_1-ta_2\bigr)^2}{\gamma^2}.
	\]
	Hence
	\[
	\Psi_{t,X}(u_3^*)
	=
	(4b_1b_2+\gamma^2)
	\left[
	\rho^2
	+
	\frac{\bigl((1-t)b_1-tb_2\bigr)^2}{4b_1b_2}
	-
	\frac{\bigl(X-(1-t)a_1-ta_2\bigr)^2}{\gamma^2}
	\right].
	\]
	Using \(\rho^2=t(1-t)\), we have
	\[
	4b_1b_2\rho^2+
	\bigl((1-t)b_1-tb_2\bigr)^2
	=
	\bigl((1-t)b_1+tb_2\bigr)^2.
	\]
	Therefore,
	\[
	\Psi_{t,X}(u_3^*)
	=
	(4b_1b_2+\gamma^2)
	\left[
	\frac{\bigl((1-t)b_1+tb_2\bigr)^2}{4b_1b_2}
	-
	\frac{\bigl(X-(1-t)a_1-ta_2\bigr)^2}{\gamma^2}
	\right].
	\]
	For fixed \(X\), define
	\[
	\Phi_X(t)
	=
	\frac{\bigl((1-t)b_1+tb_2\bigr)^2}{4b_1b_2}
	-
	\frac{\bigl(X-(1-t)a_1-ta_2\bigr)^2}{\gamma^2}.
	\]
	Then
	\[
	\Psi_{t,X}(u_3^*)
	=
	(4b_1b_2+\gamma^2)\Phi_X(t).
	\]
	Set
	\[
	\Delta=a_2-a_1,
	\qquad
	d=b_1-b_2.
	\]
	Then
	\[
	\Phi_X(t)
	=
	\frac{(b_1-dt)^2}{4b_1b_2}
	-
	\frac{(X-a_1-\Delta t)^2}{\gamma^2}.
	\]
	Define
	\[
	\Lambda
	=
	4b_1b_2\Delta^2-d^2\gamma^2=
	4b_1b_2(a_2-a_1)^2-(b_1-b_2)^2\gamma^2.
	\tag{*}
	\]
	Then
	\[
	\Phi_X''(t)
	=
	-\frac{\Lambda}{2b_1b_2\gamma^2}.
	\]
	Hence \(\Phi_X\) is strictly concave when \(\Lambda>0\), affine when
	\(\Lambda=0\), and strictly convex when \(\Lambda<0\).
	For fixed \(X\), the interior critical point is admissible precisely
	when
	\(
	|u_3^*|<S.
	\)
	Equivalently,
	\(
	H_X(t)>0,
	\)
	where
	\[
	H_X(t)
	=
	t(1-t)
	-\frac{(X-a_1-\Delta t)^2}{\gamma^2}
	-\frac{\gamma^2\bigl((1-t)b_1-tb_2\bigr)^2}
	{16b_1^2b_2^2}.
	\]
	Since \(H_X\) is strictly concave, its positive set, when nonempty,
	is an interval
	\(
	I_X=(\alpha_X,\beta_X),
	\)
	whose endpoints satisfy
	\(
	|u_3^*|=S.
	\)
	
	Suppose first that \(\Lambda=0\). Then \(\Phi_X\) is affine in \(t\).
	Hence its maximum over the closure of the admissible interval
	\([\alpha_X,\beta_X]\) is attained at an endpoint, unless \(\Phi_X\)
	is constant. In the latter case, the interior values coincide with
	the endpoint values. Since
	\(
	|u_3^*|=S
	\)
	at
	\(
	t=\alpha_X,\beta_X,
	\)
	the interior critical branch produces no new exposed boundary
	component.
	
	If \(\Lambda<0\), then \(\Phi_X\) is strictly convex. Therefore its
	maximum over \([\alpha_X,\beta_X]\) is attained at one of the two
	endpoints. Again, \(|u_3^*|=S\) there, so the maximizing point reduces
	to one of the endpoint branches. Hence no additional upper-boundary
	component is produced.
	
		It remains to consider \(\Lambda>0\). In this case \(\Phi_X\) is
	strictly concave. Hence, if an admissible stationary point lies in the
	interior of the admissible \(t\)-interval, it is the unique maximizer
	of \(\Phi_X\) over that interval and gives a value strictly larger than
	those at the endpoints. We now determine when such an admissible
	interior stationary point exists. It must satisfy
	\[
	0=\Phi_X'(t)
	=
	-\frac{d(b_1-dt)}{2b_1b_2}
	+
	\frac{2\Delta(X-a_1-\Delta t)}{\gamma^2},
	\]
	or equivalently,
	\begin{align}\label{adeq}
		 4b_1b_2\Delta(X-a_1-\Delta t)
		 =
		 \gamma^2d(b_1-dt).
	\end{align}

	Along the stationary branch, substituting (\ref{adeq}) into \(H_X(t)\) gives
	\[
	16b_1^2b_2^2H_X(t)=G(t),
	\]
	where
	\[
    \begin{aligned}
    	G(t)
    	={}&16b_1^2b_2^2t(1-t)
    	-\gamma^2\bigl((1-t)b_1-tb_2\bigr)^2
    	-\frac{\gamma^2(b_1-b_2)^2}{(a_2-a_1)^2}
    	\bigl((1-t)b_1+tb_2\bigr)^2.
    \end{aligned}
	\]
	Thus the stationary point is admissible precisely when
	\(
	G(t)>0.
	\)

	Let
	\(
	t_0=\frac{b_1}{b_1+b_2}.
	\)
	Then
	\[
	G(t_0)
	=
	\frac{4b_1^2b_2^2}
	{(b_1+b_2)^2\Delta^2}\Lambda>0,
	\qquad
	G(0)<0,
	\qquad
	G(1)<0.
	\]
	Since \(G\) is continuous and strictly concave, and
   \[
    G(0)<0,\qquad G(t_0)>0,\qquad G(1)<0,
   \]
   there exist unique
   \(
    t_-\in(0,t_0),~
   t_+\in(t_0,1)
   \)
   such that
   \(
   G(t_-)=G(t_+)=0.
   \)
   Moreover,
   \[
   G(t)>0
   \quad\Longleftrightarrow\quad
   t\in(t_-,t_+).
   \]
	At \(t=t_\pm\), one has
	\(
	|u_3^*|=S.
	\)
	
	Moreover, since
	\(
	t_-<t_0=\frac{b_1}{b_1+b_2}<t_+,
	\)
	we have
	\(
	A_{t_-}-B_{t_-}>0,~
	A_{t_+}-B_{t_+}<0.
	\)
	Hence
	\(
	u_3^*(t_-)=S,~
	u_3^*(t_+)=-S.
	\)
	Thus the interior critical branch joins the \(u_3=S\) endpoint branch
	at \(t=t_-\) and the \(u_3=-S\) endpoint branch at \(t=t_+\).
	
	At the stationary point, we have
	\[
	X-a_1-\Delta t
	=
	\frac{\gamma^2d(b_1-dt)}
	{4b_1b_2\Delta}.
	\]
	Substituting this into \(\Phi_X(t)\) gives
	\[
	\Phi_X(t)
	=
	\frac{\Lambda(b_1-dt)^2}
	{16b_1^2b_2^2\Delta^2}.
	\]
	Moreover,
	\[
	dX+b_2a_1-b_1a_2
	=
	-\frac{\Lambda}{4b_1b_2\Delta}(b_1-dt).
	\]
	Hence
	\[
	\Phi_X(t)
	=
	\frac{(dX+b_2a_1-b_1a_2)^2}{\Lambda},
	\]
	and therefore
	\[
	Y^2
	=
	\frac{4b_1b_2+\gamma^2}{\Lambda}
	(dX+b_2a_1-b_1a_2)^2.
	\]
	Since \(b_1-dt=(1-t)b_1+tb_2>0\), the quantity
	\(dX+b_2a_1-b_1a_2\) has a fixed sign along the admissible stationary
	branch. Thus this branch is a straight line segment. More explicitly, along this branch,
	\[
	X(t)=a_1+\Delta t+
	\frac{\gamma^2d(b_1-dt)}{4b_1b_2\Delta},
	\qquad
	X'(t)=\frac{\Lambda}{4b_1b_2\Delta}\ne0.
	\]
	Thus \(t\in[t_-,t_+]\) traces the entire segment exactly once,
	including its two joining points.
	
	If \(a_1=a_2\), then
	\[
	\Lambda=-(b_1-b_2)^2\gamma^2\le0,
	\]
	so this case is already covered by the preceding analysis. If, in addition, \(b_1=b_2\), then \(\Lambda=0\) and the two
	auxiliary ellipses coincide, \(E_+=E_-\). Hence no additional
	tangent segment occurs.

   For the following common-tangent calculation, assume \(\Lambda>0\). Set
   \(
   s=b_1+b_2.
   \)
   Then
   \(
   s^2-d^2=4b_1b_2.
   \)
   The stationary branch obtained above is contained in the line
   \[
   Y=q(dX+b_2a_1-b_1a_2),
   \qquad
   q^2=\frac{4b_1b_2+\gamma^2}{\Lambda}.
   \]
   Here the sign of \(q\) is chosen so that
   \(q(dX+b_2a_1-b_1a_2)\geq0\) along the admissible branch.

	We now show that this line is a common tangent of the two auxiliary
	ellipses.
	
	The ellipse \(E_+\) has center
	\[
	C_+
	=
	\left(
	\frac{a_1+a_2}{2},
	\frac{d}{2}
	\right)
	\]
	and shape matrix
	\[
	Q_+
	=
	\frac14\left(v_+v_+^T+\gamma^2I\right),
	\qquad
	v_+=(\Delta,-s).
	\]
    Thus the boundary of \(E_+\) may be written as
	\[
	(z-C_+)^TQ_+^{-1}(z-C_+)=1.
	\]
	Writing the candidate line in the form
	\[
	n^Tz=r,
	\]
	where
	\(
	n=(-qd,1),
    ~
	r=q(b_2a_1-b_1a_2).
	\)
	A line \(n^Tz=r\) is tangent to \(E_+\) precisely when
	\[
	|r-n^TC_+|^2=n^TQ_+n.
	\]
	Set
	\[
	c=\frac{a_1+a_2}{2}.
	\]
	Since
	\[
	dc+b_2a_1-b_1a_2
	=
	-\frac{s\Delta}{2},
	\]
	we obtain
	\[
	r-n^TC_+
	=
	-\frac{qs\Delta+d}{2}.
	\]
	Hence
	\[
	4|r-n^TC_+|^2
	=
	(qs\Delta+d)^2.
	\]
	
	On the other hand,
	\[
	n^Tv_+=-(qd\Delta+s),
	\qquad
	\|n\|^2=q^2d^2+1,
	\]
	and therefore
	\[
	4n^TQ_+n
	=
	(qd\Delta+s)^2
	+
	\gamma^2(q^2d^2+1).
	\]
	The difference between the two sides is
	\[
	\begin{aligned}
		&~\quad(qs\Delta+d)^2
		-
		(qd\Delta+s)^2
		-
		\gamma^2(q^2d^2+1)\\
		&=
		q^2\bigl((s^2-d^2)\Delta^2-d^2\gamma^2\bigr)
		-(s^2-d^2+\gamma^2)\\
		&=
		q^2\Lambda-(4b_1b_2+\gamma^2)
		=0.
	\end{aligned}
	\]
	Thus
	\[
	|r-n^TC_+|^2=n^TQ_+n,
	\]
	and the line is tangent to \(E_+\).
	
	Similarly, \(E_-\) has center
	\[
	C_-
	=
	\left(
	\frac{a_1+a_2}{2},
	-\frac{d}{2}
	\right)
	\]
	and shape matrix
	\[
	Q_-
	=
	\frac14\left(v_-v_-^T+\gamma^2I\right),
	\qquad
	v_-=(\Delta,s).
	\]
	In this case,
	\[
	r-n^TC_-
	=
	\frac{d-qs\Delta}{2},
	\]
	while
	\[
	n^Tv_-=s-qd\Delta.
	\]
	Hence
	\begin{align*}
		&~\quad4|r-n^TC_-|^2-4n^TQ_-n\\
		&=
		(d-qs\Delta)^2
		-(s-qd\Delta)^2
		-\gamma^2(q^2d^2+1)\\
		&=
		q^2\Lambda-(4b_1b_2+\gamma^2)
		=0.
    \end{align*}
	Therefore the same line is tangent to \(E_-\).
	
	Finally, we have already shown that
	\(
	u_3^*(t_-)=S,
	~
	u_3^*(t_+)=-S.
	\)
	Hence the two endpoints of the admissible stationary segment coincide
	with the corresponding endpoint branches on \(E_+\) and \(E_-\),
	respectively. Therefore the stationary branch is precisely the common
	tangent segment joining the two relevant ellipse arcs.
	
	For fixed \(t,X\), put \(c=A_t-B_t\) and \(R=\gamma\rho S\).
	The larger of the two attainable endpoint heights is
	\[
	\max\{|c+R|,|c-R|\}
	=|c|+R
	=\max\{c+R,-c+R\}.
	\]
	
   The last two terms are the upper heights of the two auxiliary disks
   at the same \(t,X\). Maximizing them over all feasible \(t\)
   therefore gives the larger of the two upper elliptic arcs.
   
   It remains to compare these endpoint heights with all admissible
   interior critical values. For fixed \(X\), the squared interior
   critical heights are
   \((4b_1b_2+\gamma^2)\Phi_X(t)\) on the interval where \(H_X(t)>0\).
   If this interval is empty, there is no interior contribution.
   Otherwise, if \(\Lambda\le0\), their maximum over its closure
   occurs at an endpoint, where \(|u_3^*|=S\), so it is already
   an endpoint height.
   
   If \(\Lambda>0\), strict concavity shows that the maximum occurs
   either at an admissible stationary point or at an endpoint of that
   interval. In the latter case it again reduces to an endpoint height.
   In the former case, the stationary value is the maximum of the
   quadratic \(\Phi_X\) over all real \(t\), and its point lies on
   the common tangent segment identified above.
   Its comparison with all endpoint heights follows from the
   upper-supporting property established next.

	When \(\Lambda>0\), the common tangent found above is the upper
	supporting line of both auxiliary ellipses. Indeed, the sign
	choice for \(q\) gives \(q\Delta<0\), while
	\(q^2=(4b_1b_2+\gamma^2)/\Lambda\) and
	\(s^2-d^2=4b_1b_2\) give \(|q|s|\Delta|>|d|\).
	Thus \(r-n^TC_+>0\) and \(r-n^TC_->0\); together with the
	tangency equalities above, this places both ellipses in
	\(n^Tz\leq r\).
	
	Consequently, the upper boundary of $\bild^{+}(A)$ is the upper boundary of $\operatorname{conv}(E_+\cup E_-)$, and is therefore concave. In fact, let $m(X)$ and $h(X)$ denote the midpoint and the half-length of the vertical section of $E_+$ at $X$. Since $E_-$ is the reflection of $E_+$, the larger of the two upper arcs is $h(X)+|m(X)|$, where $h$ is concave and $m$ is affine. The extreme points of $E_+$ in the $X$-direction have heights
	\[
	\frac{b_1-b_2}{2}\mp\frac{(a_2-a_1)(b_1+b_2)}{2\sqrt{(a_2-a_1)^2+\gamma^2}},
	\]
	which have opposite signs if and only if $\Lambda>0$. Thus, if $\Lambda\le0$, then $m$ does not change sign and one upper arc dominates the other; if $\Lambda>0$, the common tangent segment bridges the two arcs.
	
    For the lower boundary, keep \(c=A_t-B_t\) and \(R=\gamma\rho S\)
    for fixed feasible \(t,X\). The endpoint calculation gives the
    lowest attainable height
    \[
    \min\{|c+R|,|c-R|\}=\bigl||c|-R\bigr|.
    \]
    The vertical section of \(E_+\) at \(X\) has endpoints
    \(m(X)-h(X)\) and \(m(X)+h(X)\). If this section lies entirely
    above or below the real axis, its distance from the axis is the
    actual lowest height. If it meets the real axis, some feasible
    \(t\) satisfies \(|c|\leq R\). At the endpoints of the feasible
    \(t\)-interval, \(R=0\); continuity therefore gives a feasible
    \(t\) with \(|c|=R\), so height zero is actually attained.
    Consequently, the lower-boundary function is
	\[
	\ell(X)=\max\{0,|m(X)|-h(X)\}.
	\]
	Since \(m\) is affine and \(h\) is concave, this function is convex. Together with Lemma~\ref{connected}, this proves that $\bild^{+}(A)$ is convex.
	
	\begin{proposition}\label{prop:pure-j}
		Let
		\[
		A=
		\begin{pmatrix}
			a_1+b_1i&\gamma j\\
			0&a_2+b_2i
		\end{pmatrix},
		\qquad a_{1},a_{2}\in \mathbb{R}, b_1,b_2,\gamma>0,
		\]
		and define
		\[
		\Lambda
		=
		4b_1b_2(a_2-a_1)^2-(b_1-b_2)^2\gamma^2.
		\]
		Then \(\bild^+(A)\) is convex. More precisely, let \(E_+\) and \(E_-\)
		be the two auxiliary ellipses determined above.
		\begin{enumerate}
			\item If \(\Lambda\leq0\), the interior critical branch contributes
			no additional upper-boundary component.
			
			\item If \(\Lambda>0\), the interior critical branch is a nontrivial
			common tangent segment joining the relevant arcs of \(E_+\) and
			\(E_-\).
		\end{enumerate}
	\end{proposition}
	
\begin{remark}
	For the matrices in Proposition~\ref{prop:pure-j},
	\(W_{\H}(A)\) is convex if and only if \(a_1=a_2\) and \(b_1=b_2\).
	Necessity follows from the preceding boundary description.
	The vertical sections at both extreme real parts must meet the real axis.
	If \(a_1=a_2=a\) and \(b_1=b_2=b\), the two auxiliary ellipses
	coincide and are centered on the real axis, so the full bild is an
	elliptical disk. Hence \(W_{\H}(A)\) is convex by
	\cite[Theorem~3.1]{CDM20star}. In this case, \(A\) is also unitarily similar to a real matrix.
\end{remark}

\section{The Generic Case: A Distance-to-a-Disk Representation}
\label{Section:Generic}

In this section, we first show that the upper boundary of the upper bild
of an \(n\times n\) quaternionic matrix is concave. We then consider the
generic \(2\times2\) case and prove that the lower-boundary function is
convex by expressing it as the distance from an affinely moving point to
a disk with concave radius. Together with Lemma~\ref{connected}, these
boundary properties imply the convexity of the upper bild in the
\(2\times2\) case. Lemma~\ref{2x2} then extends the result to arbitrary
dimension.

The following fact is implicit in \cite{K51,K08} and was later used explicitly by \cite{ST1}. For completeness, we provide a direct proof.

\begin{lemma}\label{upconcave}
	Let \(A\in M_n(\mathbb H)\), and let
	\(\chi_A\in M_{2n}(\mathbb C)\) be its complex representation.
	Define the upper bild of \(A\) by
	\[
	\bild^{+}(A)
	=
	\left\{
	\operatorname{Re}(x^{*}Ax)
	+i\left|\operatorname{Im}_{\mathbb H}(x^{*}Ax)\right|
	:
	x\in S_{\H}^{n-1}
	\right\},
	\]
	and let
	\[
	C=W_{\mathbb C}(\chi_A)
	=
	\left\{
	z^{*}\chi_Az:
	z\in\mathbb C^{2n},\ z^{*}z=1
	\right\}.
	\]
	
	Let
	\[
	I_B=\{\operatorname{Re}w:w\in \bild^{+}(A)\},
	\qquad
	I_C=\{\operatorname{Re}w:w\in C\},
	\]
	and define the corresponding upper-boundary functions by
	\[
	u_B(t)=\max\{y\geq 0:t+iy\in \bild^{+}(A)\},
	\qquad t\in I_B,
	\]
	and
	\[
	u_C(t)=\max\{y:t+iy\in C\},
	\qquad t\in I_C.
	\]
	
	Then
	\[
	\bild^{+}(A)\subseteq C,
	\qquad
	I_B=I_C,
	\qquad
	u_B(t)=u_C(t)
	\quad\text{for every }t\in I_B=I_C.
	\]
	Consequently, \(u_B\) is a concave function. Equivalently, the
	upper boundary of \(\bild^{+}(A)\) is concave.
\end{lemma}

\begin{proof}
	We use the complex representation convention
	\[
	q=m_1+m_2j
	\quad\longmapsto\quad
	\chi_q=
	\begin{pmatrix}
		m_1&m_2\\
		-\overline{m_2}&\overline{m_1}
	\end{pmatrix},
	\qquad m_1,m_2\in\mathbb C.
	\]
	The two eigenvalues of \(\chi_q\) are
	\(
	\operatorname{Re}(q)
	\pm i\left|\operatorname{Im}_{\mathbb H}(q)\right|,
	\)
	where
	\(
	\left|\operatorname{Im}_{\mathbb H}(q)\right|
	=
	\sqrt{(\operatorname{Im}m_1)^2+|m_2|^2}.
	\)

	First, we prove that \(\bild^{+}(A)\subseteq C\). Let \(x\in\mathbb H^n\) be a quaternionic unit vector and set
	\(
	q=x^{*}Ax.
	\)
	Let \(\chi_x\in\mathbb C^{2n\times 2}\) be the complex representation
	of \(x\). Since \(x^{*}x=1\), we have
	\(
	\chi_x^{*}\chi_x=I_2.
	\)
	Moreover,
	\(
	\chi_x^{*}\chi_A\chi_x
	=
	\chi_{x^{*}Ax}
	=
	\chi_q.
	\)
	The upper-bild representative
	\(
	\widehat q
	=
	\operatorname{Re}(q)
	+i\left|\operatorname{Im}_{\mathbb H}(q)\right|
	\)
	is one of the two eigenvalues of \(\chi_q\). Hence there exists a unit
	vector \(v\in\mathbb C^2\) such that
	\(
	v^{*}\chi_qv=\widehat q.
	\)
	Because \(\chi_x^{*}\chi_x=I_2\), the vector \(\chi_xv\) is a unit
	vector in \(\mathbb C^{2n}\). Therefore
	\[
	\widehat q
	=
	v^{*}\chi_x^{*}\chi_A\chi_xv
	=
	(\chi_xv)^{*}\chi_A(\chi_xv)
	\in W_{\mathbb C}(\chi_A).
	\]
	Thus
	\(
	\bild^{+}(A)\subseteq C.
	\)
	In particular,
	\(
	I_B\subseteq I_C
	\)
	and
	\(
	u_B(t)\leq u_C(t),
	~ t\in I_B.
	\)
	
    Second, we prove that \(I_C\subseteq I_B\) and \(u_B(t)\geq u_C(t)\). Fix \(t\in I_C\). Since \(C\) is compact, there exists a complex unit
	vector \(z\in\mathbb C^{2n}\) such that
	\[
	z^{*}\chi_Az=t+iu_C(t).
	\]
	Every unit vector \(z\in\mathbb C^{2n}\) can be realized as the first
	column of \(\chi_x\) for some quaternionic unit vector
	\(x\in\mathbb H^n\). For this \(x\), write
	\[
	\chi_x^{*}\chi_A\chi_x
	=
	\begin{pmatrix}
		m_1&m_2\\
		-\overline{m_2}&\overline{m_1}
	\end{pmatrix}.
	\]
	Since \(z\) is the first column of \(\chi_x\),
	\(
	m_1=z^{*}\chi_Az=t+iu_C(t).
	\)
	The corresponding quaternionic Rayleigh quotient is
	\(
	x^{*}Ax=m_1+m_2j.
	\)
	Its upper-bild coordinates are therefore
	\[
	X=\operatorname{Re}(m_1)=t,\quad
	Y
	=
	\sqrt{(\operatorname{Im}m_1)^2+|m_2|^2}
	=
	\sqrt{u_C(t)^2+|m_2|^2}.
	\]
    Let
    \(
    J_n=\begin{pmatrix}0&I_n\\-I_n&0\end{pmatrix}.
    \)
    The symplectic form of \(\chi_A\) gives
    \(J_n^*\chi_AJ_n=\overline{\chi_A}\). Consequently,
    \[
    C=W_{\mathbb C}(\chi_A)
    =W_{\mathbb C}(\overline{\chi_A})
    =\overline{W_{\mathbb C}(\chi_A)}.
    \]
    Thus \(C\) is symmetric about the real axis, and \(u_C(t)\geq0\).
    Thus \(\bild^{+}(A)\) contains a point with horizontal coordinate \(t\)
	and vertical coordinate at least \(u_C(t)\). Consequently,
	\(
	t\in I_B
	\)
	and
	\(
	u_B(t)\geq u_C(t).
	\)
	Therefore
	\(
	I_C\subseteq I_B.
	\)
	
	Combining the first step and the second step, we have 
	\(
	I_B=I_C
	\)
	and
	\(
	u_B(t)=u_C(t)\) for every
    \( t\in I_B.
	\)
	By the Toeplitz--Hausdorff theorem,
	\(
	C=W_{\mathbb C}(\chi_A)
	\)
	is convex. Hence \(u_C\), and therefore \(u_B\), is concave.
\end{proof}

We now turn to the lower boundary for \(2\times2\) quaternionic matrices.

Let \(A\in M_{2}(\mathbb H)\), and write
\[
A=H+K,
\qquad
H=\frac{A+A^{*}}{2},
\qquad
K=\frac{A-A^{*}}{2},
\]
so that \(H^{*}=H\) and \(K^{*}=-K\).  For a quaternionic unit vector
\(x\in\mathbb H^{2}\), the real and imaginary parts of the Rayleigh
quotient \(x^{*}Ax\) are determined by
\[
\operatorname{Re}(x^{*}Ax)=x^{*}Hx,
\qquad
\operatorname{Im}(x^{*}Ax)=x^{*}Kx.
\]
Consequently, after passing to the upper representative of the
quaternionic similarity class, the corresponding point in the upper
bild is
\(
\left(x^{*}Hx,\; |x^{*}Kx|\right).
\)

To reduce the number of parameters, we use the following normalization lemma.

\begin{lemma}
	\label{lem:matrix-normalization}
	Let
	\[
	H=
	\begin{pmatrix}
		h_1&0\\
		0&h_2
	\end{pmatrix},
	\qquad
	~
	K=
	\begin{pmatrix}
		\alpha&z\\
		-z^*&\beta
	\end{pmatrix},
	\]
	where \(h_1,h_2\in\mathbb R\),
	\(\alpha,\beta\in\operatorname{Im}\mathbb H\), and \(z\neq0\). Then there exists a diagonal quaternionic
	unitary matrix
	\[
	D=
	\begin{pmatrix}
		r_1&0\\
		0&r_2
	\end{pmatrix},
	\qquad
	r_1,r_2\in\operatorname{Sp}(1),
	\]
	such that
	\[
	D^*HD=H, \quad
	D^*KD
	=
	\begin{pmatrix}
		ai&\rho\\
		-\rho&b(ci+sj)
	\end{pmatrix},
	\]
	where
	\(
	\rho=|z|>0,~
	a=|\alpha|,~
	b=|\beta|,
	\)
	and
	\(
	c\in[-1,1],~
	s\geq 0,~
	c^2+s^2=1.
	\)
\end{lemma}
A proof of Lemma~\ref{lem:matrix-normalization} is given in the Supplementary Material (Lemma~2).

After a quaternionic unitary similarity and an interchange of the two
basis vectors if necessary, we may write
\[
H=
\begin{pmatrix}
	h_1&0\\
	0&h_2
\end{pmatrix},
\qquad h_1\leq h_2.
\]
If \(h_1=h_2\), then \(H=h_1I_2\), and the upper bild is a vertical
line segment, hence convex. We therefore assume \(h_1<h_2\).

By Lemma~\ref{lem:matrix-normalization}, we may also assume that \(K\)
has the normal form given above. We first consider the generic case
\(a,b,\rho,s>0\). The remaining cases will be treated later.

With this normal form, we now rewrite \(x^*Kx\) in real coordinates.

\begin{lemma}\label{A0B0}
	\label{lem:block-decomposition}
	Let
	\[
	H=
	\begin{pmatrix}
		h_1&0\\
		0&h_2
	\end{pmatrix},
	\qquad
	K=
	\begin{pmatrix}
		ai&\rho\\
		-\rho&b(ci+sj)
	\end{pmatrix},
	\]
	where
	\(
	h_1,h_2\in\mathbb R,
	~
	a,b,\rho\geq 0,
	~
	c,s\in\mathbb R,
	~
	c^2+s^2=1.
	\)
	For $0\leq\lambda\leq1$, let
	\[
	x=
	\begin{pmatrix}
		\sqrt{1-\lambda}\\
		\sqrt{\lambda}\,u
	\end{pmatrix},
	\qquad
	u\in\operatorname{Sp}(1).
	\]
	Write
	\(
	u=u_0+u_1i+u_2j+u_3k,
	~
	u_0^2+u_1^2+u_2^2+u_3^2=1,
	\)
	and define
	\[
	W=
	\begin{pmatrix}
		u_0\\
		u_3
	\end{pmatrix},
	\qquad
	V=
	\begin{pmatrix}
		u_1\\
		u_2
	\end{pmatrix},
    \qquad
	\omega_\lambda=
	\begin{pmatrix}
		(1-\lambda)a+\lambda bc\\
		\lambda bs
	\end{pmatrix},
	\]
	together with
	\[
	A_0=
	\begin{pmatrix}
		0&2abs\\
		2abs&4(\rho^2-abc)
	\end{pmatrix},
	\quad
	B_0=
	\begin{pmatrix}
		4\rho^2&2abs\\
		2abs&4(\rho^2-abc)
	\end{pmatrix}.
	\]
	Then
	\[
	x^*Hx=(1-\lambda)h_1+\lambda h_2
	\]
	and
	\[
		\begin{aligned}
			|x^*Kx|^2
			={}&
			|\omega_\lambda|^2
			+\lambda(1-\lambda)W^TA_0W+\lambda(1-\lambda)V^TB_0V+4\rho\sqrt{\lambda(1-\lambda)}
			\,\omega_\lambda^TV.
		\end{aligned}
	\]
	Moreover,
	\[
	B_0
	=
	A_0+
	\begin{pmatrix}
		4\rho^2&0\\
		0&0
	\end{pmatrix}.
	\]
\end{lemma}

\begin{proof}
	By Lemma \ref{xp}, it suffices to consider
	\[
	x=
	\begin{pmatrix}
		\sqrt{1-\lambda}\\
		\sqrt{\lambda}\,u
	\end{pmatrix}.
	\]
    By direct calculation, we have that 
	\[
	\begin{aligned}
		x^*Hx=(1-\lambda)h_1+\lambda h_2,
	\end{aligned}
	\]
	\[
	\begin{aligned}
		x^*Kx &=
		(1-\lambda)ai
		+\rho\sqrt{\lambda(1-\lambda)}\,u
		-\rho\sqrt{\lambda(1-\lambda)}\,u^*
		+\lambda u^*b(ci+sj)u\\
		&=(1-\lambda)ai
		+\lambda u^*b(ci+sj)u
		+2\rho\sqrt{\lambda(1-\lambda)}
		\,\operatorname{Im}u.
	\end{aligned}	
	\]
	For convenience, set
	\(
	\tau_\lambda=\sqrt{\lambda(1-\lambda)},
	~
	\alpha=ai,
	~
	\beta=b(ci+sj).
	\)
	Then
	\[
	x^*Kx
	=
	(1-\lambda)\alpha
	+\lambda u^*\beta u
	+2\rho\tau_\lambda\operatorname{Im}u.
	\]
    Thus,
	\[
	\begin{aligned}
		|x^*Kx|^2
		&=(1-\lambda)^2|\alpha|^2+\lambda^2|u^*\beta u|^2
		+4\rho^2\tau_\lambda^2|\operatorname{Im}u|^2\\
		&\quad+2\lambda(1-\lambda)\langle\alpha,u^*\beta u\rangle
		+4\rho\tau_\lambda(1-\lambda)\langle\alpha,\operatorname{Im}u\rangle
		+4\rho\tau_\lambda\lambda\langle u^*\beta u,\operatorname{Im}u\rangle\\
		&=(1-\lambda)^2a^2+\lambda^2b^2
		+4\rho^2\lambda(1-\lambda)(u_1^2+u_2^2+u_3^2)\\
		&\quad+2\lambda(1-\lambda)\langle\alpha,u^*\beta u\rangle
		+4\rho\tau_\lambda(1-\lambda)\langle\alpha,\operatorname{Im}u\rangle
		+4\rho\tau_\lambda\lambda\langle u^*\beta u,\operatorname{Im}u\rangle,
	\end{aligned}
	\]
	where
	\(
	|\alpha|=a,
	|u^*\beta u|=|\beta|=b,
	\)
	and
	\(
	|\operatorname{Im}u|^2=u_1^2+u_2^2+u_3^2.
	\)
	A direct quaternionic multiplication gives
	\[
	\begin{aligned}
		u^*\beta u
		={}&b\Big[
		\bigl(
		c(1-2u_2^2-2u_3^2)
		+2s(u_1u_2+u_0u_3)
		\bigr)i
		+
		\bigl(
		2c(u_1u_2-u_0u_3)
		+s(1-2u_1^2-2u_3^2)
		\bigr)j\\
		&\quad\quad\quad\quad\quad\quad\quad\quad\quad\quad\quad\quad\quad\quad\quad\quad\quad+
		\bigl(
		2c(u_1u_3+u_0u_2)
		+2s(u_2u_3-u_0u_1)
		\bigr)k
		\Big].
	\end{aligned}
	\]
	Since $\alpha=ai$, we have
	\(
		\langle\alpha,u^*\beta u\rangle
		=
		ab\Big[
		c(1-2u_2^2-2u_3^2)
		+2s(u_1u_2+u_0u_3)
		\Big].
	\)
	It follows that
	\[
	\begin{aligned}
		2\lambda(1-\lambda)
		\langle\alpha,u^*\beta u\rangle=
		2\lambda(1-\lambda)abc
		-4\lambda(1-\lambda)abc(u_2^2+u_3^2)
		+4\lambda(1-\lambda)abs(u_1u_2+u_0u_3).
	\end{aligned}
	\]
	Since
	\(
	\alpha=ai
	\)
	and
	\(
	\operatorname{Im}u=u_1i+u_2j+u_3k,
	\)
	we have
	\(
	\langle\alpha,\operatorname{Im}u\rangle=au_1
	\)
	and
    \(
    \langle u^*\beta u,\operatorname{Im}u\rangle
    =
    \langle\beta,\operatorname{Im}u\rangle=
    b(cu_1+su_2).
    \)
    Here we used the fact that \(u\) commutes with \(\operatorname{Im}u\).
    It follows that
    \begin{align*}
    	&~\quad 4\rho\tau_\lambda(1-\lambda)
    	\langle\alpha,\operatorname{Im}u\rangle
    	+4\rho\tau_\lambda\lambda
    	\langle u^*\beta u,\operatorname{Im}u\rangle
        \\
    	&=
    	4\rho\tau_\lambda(1-\lambda)au_1+4\rho\tau_\lambda\lambda b(cu_1+su_2)\\
    	&=
    	4\rho\tau_\lambda
    	\left[
    	\bigl((1-\lambda)a+\lambda bc\bigr)u_1
    	+\lambda bs\,u_2
    	\right]
    	\\
    	&=
    	4\rho\tau_\lambda\,\omega_\lambda^TV.
    \end{align*}

    Also,
	
	\begin{align}\label{olambda2}
		|\omega_\lambda|^2=(1-\lambda)^2a^2
		+\lambda^2b^2
		+2\lambda(1-\lambda)abc,
	\end{align}
	\[
	\begin{aligned}
		W^TA_0W
		&=
		\begin{pmatrix}u_0&u_3\end{pmatrix}
		\begin{pmatrix}
			0&2abs\\
			2abs&4(\rho^2-abc)
		\end{pmatrix}
		\begin{pmatrix}
			u_0\\u_3
		\end{pmatrix}=
		4abs\,u_0u_3
		+4(\rho^2-abc)u_3^2,
	\end{aligned}
	\]
	and
	\[
	\begin{aligned}
		V^TB_0V
		&=
		\begin{pmatrix}u_1&u_2\end{pmatrix}
		\begin{pmatrix}
			4\rho^2&2abs\\
			2abs&4(\rho^2-abc)
		\end{pmatrix}
		\begin{pmatrix}
			u_1\\u_2
		\end{pmatrix}=
		4\rho^2u_1^2
		+4abs\,u_1u_2
		+4(\rho^2-abc)u_2^2.
	\end{aligned}
	\]

 Combining these identities gives
	\[
	\begin{aligned}
		|x^*Kx|^2
		=
		|\omega_\lambda|^2
		+\lambda(1-\lambda)W^TA_0W
		+\lambda(1-\lambda)V^TB_0V
		+4\rho\sqrt{\lambda(1-\lambda)}
		\,\omega_\lambda^TV.
	\end{aligned}
	\]
\end{proof}

By Lemma~\ref{A0B0}, the interval of possible real parts is
\([h_1,h_2]\). We parametrize it by
\[
t(\lambda)=(1-\lambda)h_1+\lambda h_2,
\qquad 0\leq\lambda\leq1.
\]
The lower-boundary function is therefore
\begin{align}\label{minxKx}
	\ell(\lambda)
	=
	\min_{u\in\operatorname{Sp}(1)}
	|x^*Kx|.
\end{align}
Thus the problem is reduced to proving that \(\ell(\lambda)\) is convex
on \([0,1]\).

To reduce this minimization, we need the following Lemma related to positivity. 

\begin{lemma}
	\label{lem:C-positive-definite}
	Let \(A_0\) and \(B_0\) be as defined in Lemma~\ref{A0B0}.
	Let $\gamma$ denote the smallest eigenvalue of $A_0$, and define
	\(
	C=B_0-\gamma I_2.
	\)
	Then
	\[
	\gamma
	=
	2\left(
	\rho^2-abc
	-
	\sqrt{(\rho^2-abc)^2+a^2b^2s^2}
	\right).
	\]
	Moreover, $C$ is positive semi-definite. 
	Under the generic assumptions
	\(
	a>0,~b>0,~\rho>0,~s>0,
	\)
	the matrix $C$ is positive definite.

\end{lemma}

\begin{proof}
	We first compute the smallest eigenvalue of $A_0$. Consider a general
	real symmetric matrix of the form
	\[
	M=
	\begin{pmatrix}
		0&k\\
		k&d
	\end{pmatrix}.
	\]
	Its characteristic polynomial is
	\[
	\begin{aligned}
		\det(M-\mu I_2)
		&=
		\det
		\begin{pmatrix}
			-\mu&k\\
			k&d-\mu
		\end{pmatrix}=
		(-\mu)(d-\mu)-k^2=
		\mu^2-d\mu-k^2.
	\end{aligned}
	\]
	Thus its two eigenvalues are
	\[
	\mu_{\pm}
	=
	\frac{d\pm\sqrt{d^2+4k^2}}{2}.
	\]
	For the matrix $A_0$, we have
	\(
	k=2abs
	\)
	and
	\(
	d=4(\rho^2-abc).
	\)
	Therefore, the smallest eigenvalue of $A_0$ is
	\[
	\begin{aligned}
		\gamma
		&=
		\frac{
			4(\rho^2-abc)
			-
			\sqrt{
				16(\rho^2-abc)^2+16a^2b^2s^2
			}
		}{2}=
		2\left(
		\rho^2-abc
		-
		\sqrt{
			(\rho^2-abc)^2+a^2b^2s^2
		}
		\right).
	\end{aligned}
	\]
	Since $\gamma=\lambda_{\min}(A_0)$, all eigenvalues of
	$A_0-\gamma I_2$ are nonnegative. Hence
	\(
	A_0-\gamma I_2\succeq0.
	\)
	By the definition of $B_0$,
	\[
	\begin{aligned}
		C
		&=
		B_0-\gamma I_2=
		A_0-\gamma I_2
		+
		\begin{pmatrix}
			4\rho^2&0\\
			0&0
		\end{pmatrix}.
	\end{aligned}
	\]
	Both matrices on the right-hand side are positive semidefinite.
	Therefore,
	\(
	C\succeq0.
	\)
	
	We now assume the generic conditions
	\(
	a>0, b>0, \rho>0, s>0
	\)
	and prove that $C$ is positive definite.
	
	Suppose that
	\(
	v^TCv=0.
	\)
	Since
	\(
	v^TCv
	=
	v^T(A_0-\gamma I_2)v+4\rho^2v_1^2
	\) 
	is a sum of two nonnegative terms, both terms must vanish
	\(
	v^T(A_0-\gamma I_2)v=0
	\)
	and
	\(
	4\rho^2v_1^2=0.
	\)
	Because $\rho>0$, the second equality implies
	\(
	v_1=0.
	\)
	Thus
    $v=
    \begin{pmatrix}
    	0\\
    	v_2
    \end{pmatrix}.
    $
	Since $A_0-\gamma I_2\succeq0$, the equality
	\(
	v^T(A_0-\gamma I_2)v=0
	\)
	implies
	\(
	(A_0-\gamma I_2)v=0.
	\)
	Now
	\[
	A_0-\gamma I_2
	=
	\begin{pmatrix}
		-\gamma&2abs\\
		2abs&4(\rho^2-abc)-\gamma
	\end{pmatrix}.
	\]
	Therefore,
	\[
	\begin{aligned}
		(A_0-\gamma I_2)v
		&=
		\begin{pmatrix}
			-\gamma&2abs\\
			2abs&4(\rho^2-abc)-\gamma
		\end{pmatrix}
		\begin{pmatrix}
			0\\
			v_2
		\end{pmatrix}=
		\begin{pmatrix}
			2abs\,v_2\\
			\bigl(4(\rho^2-abc)-\gamma\bigr)v_2
		\end{pmatrix}.
	\end{aligned}
	\]
	The first component must be zero, so
	\(
	2abs\,v_2=0.
	\)
	Since
	\(
	a>0, b>0, s>0,
	\)
	we conclude that
	\(
	v_2=0.
	\)
	Hence
	\(
	v=0.
	\)
	Therefore,
	\(
	C\succ0.
	\)
\end{proof}

Now, we minimize over \(W\).

\begin{lemma}
	\label{lem:elimination-W}
	Let
	\[
	q_\lambda:=\lambda(1-\lambda),
	\qquad
	\gamma:=\lambda_{\min}(A_0),
	\qquad
	C:=B_0-\gamma I_2.
	\]
	Suppose that
	\[
	\begin{aligned}
		|x^*Kx|^2
		={}&
		|\omega_\lambda|^2
		+q_\lambda W^TA_0W
		+q_\lambda V^TB_0V+
		4\rho\sqrt{q_\lambda}\,
		\omega_\lambda^TV,
	\end{aligned}
	\]
	Then the lower-boundary function satisfies
	\[
		\begin{aligned}
			\ell(\lambda)^2
			=
			|\omega_\lambda|^2+q_\lambda\gamma
			+
			\min_{|V|\leq1}
			\left(
			q_\lambda V^TCV
			+
			4\rho\sqrt{q_\lambda}\,
			\omega_\lambda^TV
			\right).
		\end{aligned}
	\]
\end{lemma}

\begin{proof}
	By (\ref{minxKx}), we have
	\[
	\ell(\lambda)^2
	=
	\min_{u\in Sp(1)}|x^*Kx|^2.
	\]
	With \(W\) and \(V\) as defined in Lemma~\ref{A0B0},
	the condition \(u\in Sp(1)\) is equivalent to
	\(
	|W|^2+|V|^2=1.
	\)
	Fix \(V\) with \(|V|\leq1\). Then \(W\) is constrained by
	\(
	|W|^2=1-|V|^2.
	\)
	Since \(\gamma=\lambda_{\min}(A_0)\), 
	\[
	\min_{|W|^2=1-|V|^2} W^TA_0W = \gamma|W|^2
	=
	\gamma(1-|V|^2).
	\]
	Minimizing first over \(W\) consequently gives
	\[
	\begin{aligned}
		\ell(\lambda)^2
		=
		\min_{|V|\leq1}
		\Bigl[
		&|\omega_\lambda|^2
		+q_\lambda\gamma(1-|V|^2)
		+q_\lambda V^TB_0V+
		4\rho\sqrt{q_\lambda}\,
		\omega_\lambda^TV
		\Bigr].
	\end{aligned}
	\]
	Since
	\[
	\begin{aligned}
		\gamma(1-|V|^2)+V^TB_0V
		&=
		\gamma-\gamma V^TV+V^TB_0V
		\\
		&=
		\gamma+V^T(B_0-\gamma I_2)V
		\\
		&=
		\gamma+V^TCV,
	\end{aligned}
	\]
	we obtain
	\[
	\begin{aligned}
		\ell(\lambda)^2
		=
		|\omega_\lambda|^2+q_\lambda\gamma
		+
		\min_{|V|\leq1}
		\left(
		q_\lambda V^TCV
		+
		4\rho\sqrt{q_\lambda}\,
		\omega_\lambda^TV
		\right).
	\end{aligned}
	\]
\end{proof}

The following identity converts this minimum into a distance.

\begin{lemma}
	\label{lem:key-identity-distance}
	Assume
	\(
	a>0, b>0, \rho>0, s>0,
	c^2+s^2=1.
	\)
	Define \(A_{0}\), \(B_{0}\), \(C\), \(\gamma\), \(q_{\lambda}\), \(\omega_{\lambda}\) as above. 
	Then
	\begin{align}\label{wlcinverse}
		|\omega_\lambda|^2+q_\lambda\gamma
		=
		4\rho^2\omega_\lambda^TC^{-1}\omega_\lambda.
   \end{align}
   Consequently, by Lemma~\ref{lem:elimination-W},
   \[
   \ell(\lambda)
   =
   \operatorname{dist}_C
   \left(
   p(\lambda),
   \sqrt{\lambda(1-\lambda)}\,D
   \right),
   \]
	where
	\(
	p(\lambda)=-2\rho C^{-1}\omega_\lambda,
	\)
	\(
	D=\{V\in\mathbb R^2:|V|\leq1\},
	\)
	and
	\(
	\|z\|_C=\sqrt{z^TCz}.
	\)
\end{lemma}

\begin{proof}
	
	Set
	\(
	R=\rho^2-abc
	\)
	and
	\(
	\Delta=\sqrt{R^2+a^2b^2s^2}.
	\)
	Then
	\(
	\gamma=2(R-\Delta).
	\)
	Define also
	\(
	\xi=\rho^2+abc+\Delta
	\)
	and
	\(
	\eta=\rho^2-abc+\Delta=R+\Delta.
	\)
	Because $a,b,s>0$, we have
	\(
	a^2b^2s^2>0.
	\)
	Hence
	\(
	\Delta
	=
	\sqrt{R^2+a^2b^2s^2}
	>
	|R|.
	\)
	It follows that
	\(
	\eta=R+\Delta>0.
	\)

	By definition,
	\(
	C=B_0-\gamma I_2,
	\)
	that is
	\[
	C
	=
	\begin{pmatrix}
		4\rho^2-\gamma&2abs\\
		2abs&4(\rho^2-abc)-\gamma
	\end{pmatrix}.
	\]
	For the $(1,1)$ entry, we calculate
	\[
	\begin{aligned}
		4\rho^2-\gamma
		=
		4\rho^2-2(R-\Delta)
		=
		4\rho^2-2\rho^2+2abc+2\Delta
		=
		2(\rho^2+abc+\Delta)
		=
		2\xi.
	\end{aligned}
	\]
	For the $(2,2)$ entry,
	\[
	\begin{aligned}
		4(\rho^2-abc)-\gamma
		=
		4R-2(R-\Delta)
		=
		2(R+\Delta)
		=
		2\eta.
	\end{aligned}
	\]
	Thus
	\[
		C
		=
		2
		\begin{pmatrix}
			\xi&abs\\
			abs&\eta
		\end{pmatrix}.
	\]
	From the preceding formula,
	\[
	\det C
	=
	4\left(\xi\eta-a^2b^2s^2\right).
	\]
	We now simplify the expression inside the parentheses.
	First,
	\[
	\begin{aligned}
		\xi\eta
		&=
		(\rho^2+abc+\Delta)
		(\rho^2-abc+\Delta)\\
		&=
		\bigl((\rho^2+\Delta)+abc\bigr)
		\bigl((\rho^2+\Delta)-abc\bigr)\\
		&=
		(\rho^2+\Delta)^2-a^2b^2c^2\\
		&=
		\rho^4+2\rho^2\Delta+\Delta^2-a^2b^2c^2.
	\end{aligned}
	\]
	Also we have that
	\[
	\begin{aligned}
	\Delta^2
	&=(\rho^2-abc)^2+a^2b^2s^2\\
	&=\rho^4-2\rho^2abc+a^2b^2c^2+a^2b^2s^2\\
	&=
	\rho^4-2\rho^2abc+a^2b^2(c^2+s^2)\\
	&=\rho^4-2\rho^2abc+a^2b^2.
	\end{aligned}
	\]
	Therefore,
	\[
	\begin{aligned}
		\xi\eta-a^2b^2s^2
		&=
		\rho^4+2\rho^2\Delta+\Delta^2
		-a^2b^2c^2-a^2b^2s^2\\
		&=
		\rho^4+2\rho^2\Delta+\Delta^2
		-a^2b^2(c^2+s^2)\\
		&=
		\rho^4+2\rho^2\Delta+\Delta^2-a^2b^2\\
		&=
		\rho^4+2\rho^2\Delta
		+\left(
		\rho^4-2\rho^2abc+a^2b^2
		\right)
		-a^2b^2\\
		&=
		2\rho^4+2\rho^2\Delta-2\rho^2abc\\
		&=
		2\rho^2(\rho^2-abc+\Delta)\\
		&=
		2\rho^2\eta.
	\end{aligned}
	\]
	Consequently,
	\[
		\det C
		=
		4(2\rho^2\eta)
		=
		8\rho^2\eta.
	\]
	Since
	\(
	\rho>0
	\)
	and
	\(
	\eta>0,
	\)
	we have
	\(
	\det C>0.
	\)
    Now we can compute $C^{-1}$ explicitly.
	The inverse formula for a $2\times2$ matrix gives
	\[
	\begin{aligned}
	C^{-1}=
	\frac{1}{\det C}
	\begin{pmatrix}
		2\eta&-2abs\\
		-2abs&2\xi
	\end{pmatrix}=
		\frac{1}{8\rho^2\eta}
		\begin{pmatrix}
			2\eta&-2abs\\
			-2abs&2\xi
		\end{pmatrix}=
			\frac{1}{4\rho^2\eta}
			\begin{pmatrix}
				\eta&-abs\\
				-abs&\xi
			\end{pmatrix}.
	\end{aligned}
	\]

	Let
	\[
	e_0=
	\begin{pmatrix}
		a\\
		0
	\end{pmatrix},
	\qquad
	e_1=
	\begin{pmatrix}
		bc\\
		bs
	\end{pmatrix}.
	\]
	Then
	\[
	\begin{aligned}
		\omega_\lambda
		=\begin{pmatrix}
			(1-\lambda)a+\lambda bc\\
			\lambda bs
		\end{pmatrix}
		=
		(1-\lambda)
		\begin{pmatrix}
			a\\
			0
		\end{pmatrix}
		+
		\lambda
		\begin{pmatrix}
			bc\\
			bs
		\end{pmatrix}
		=
		(1-\lambda)e_0+\lambda e_1.
		\end{aligned}
	\]
    Now, we want to prove (\ref{wlcinverse}).  By calculation, we have that
	\begin{align}\label{Gram1}
		4 \rho^{2} e_0^TC^{-1}e_0
		&=
		4 \rho^{2} \begin{pmatrix}
			a&0
		\end{pmatrix}
		\frac{1}{4\rho^2\eta}
		\begin{pmatrix}
			\eta&-abs\\
			-abs&\xi
		\end{pmatrix}
		\begin{pmatrix}
			a\\
			0
		\end{pmatrix}=
		4 \rho^{2}
		\begin{pmatrix}
			a&0
		\end{pmatrix}
		\frac{1}{4\rho^2\eta}
		\begin{pmatrix}
			a\eta\\
			-a^2bs
		\end{pmatrix}=
		4 \rho^{2}
		\frac{a^2\eta}{4\rho^2\eta}=
		a^{2}
	\end{align}
	\[
	\begin{aligned}
		e_0^TC^{-1}e_1
		&=
		\begin{pmatrix}
			a&0
		\end{pmatrix}
		\frac{1}{4\rho^2\eta}
		\begin{pmatrix}
			\eta&-abs\\
			-abs&\xi
		\end{pmatrix}
		\begin{pmatrix}
			bc\\
			bs
		\end{pmatrix}
		&=
		\begin{pmatrix}
			a&0
		\end{pmatrix}
		\frac{1}{4\rho^2\eta}
		\begin{pmatrix}
			\eta bc-ab^2s^2\\
			-ab^2cs+\xi bs
		\end{pmatrix}=
		\frac{abc\eta-a^2b^2s^2}
		{4\rho^2\eta}.
	\end{aligned}
	\]
	Therefore,
	\[
	4\rho^2e_0^TC^{-1}e_1
	=
	abc-\frac{a^2b^2s^2}{\eta}.
	\]
	We now simplify the last fraction. Since
	\(
	a^2b^2s^2=\Delta^2-R^2=	(\Delta-R)(\Delta+R),
	\)
	\(
	\eta=R+\Delta,
	\)
	and
	\(
	\gamma = 2 (R-\Delta)
	\),
	it follows that
	\begin{align}\label{Gram2}
		4\rho^2e_0^TC^{-1}e_1
		=
		abc+\frac{\gamma}{2}.
	\end{align}
	Since $C^{-1}$ is symmetric,
	\begin{align}\label{Gram3}
		4\rho^2e_1^TC^{-1}e_0
		=
		abc+\frac{\gamma}{2}.
	\end{align}
	By calculation, 
	\begin{align*}
		e_1^TC^{-1}e_1
		=
		\begin{pmatrix}
			bc&bs
		\end{pmatrix}
		\frac{1}{4\rho^2\eta}
		\begin{pmatrix}
			\eta&-abs\\
			-abs&\xi
		\end{pmatrix}
		\begin{pmatrix}
			bc\\
			bs
		\end{pmatrix}
		=
		\frac{b^2}{4\rho^2\eta}
		\left(
		c^2\eta+s^2\xi-2abcs^2
		\right).
    \end{align*}
	We now simplify the expression in parentheses. Since
	\(
	\xi=\rho^2+abc+\Delta
	\)
	and
	\(
	\eta=\rho^2-abc+\Delta,
	\)
	we have
	\[
	\begin{aligned}
		\xi-\eta
		&=
		(\rho^2+abc+\Delta)
		-(\rho^2-abc+\Delta)=
		2abc.
	\end{aligned}
	\]
	Therefore,
	\[
	\begin{aligned}
		c^2\eta+s^2\xi-2abcs^2
		=
		c^2\eta+s^2(\eta+2abc)-2abcs^2
		=
		c^2\eta+s^2\eta
		=
		(c^2+s^2)\eta
		=
		\eta.
	\end{aligned}
	\]
	Consequently,
	\begin{align}\label{Gram4}
		4 \rho^{2} e_1^TC^{-1}e_1
		=
		4 \rho^{2} \frac{b^2\eta}{4\rho^2\eta}=
		b^{2}.
	\end{align}
	Since
	\[
	\omega_\lambda=(1-\lambda)e_0+\lambda e_1,
	\]
	by (\ref{Gram1})--(\ref{Gram4}), we have that 
	\begin{align*}
		4\rho^2\omega_\lambda^TC^{-1}\omega_\lambda
		&=
		4\rho^2 \bigl((1-\lambda)e_0+\lambda e_1\bigr)^T
		C^{-1}
		\bigl((1-\lambda)e_0+\lambda e_1\bigr)\\
		&=4\rho^2(1-\lambda)^2e_0^TC^{-1}e_0+4\rho^2
		\lambda(1-\lambda)e_0^TC^{-1}e_1+
		4\rho^2\lambda(1-\lambda)e_1^TC^{-1}e_0+4\rho^2
		\lambda^2e_1^TC^{-1}e_1\\
		&=4\rho^2(1-\lambda)^2e_0^TC^{-1}e_0+4\rho^2
		\cdot 2\lambda(1-\lambda)e_0^TC^{-1}e_1+
		4\rho^2 \lambda^2e_1^TC^{-1}e_1\\
		&=
		(1-\lambda)^2a^2+
		2\lambda(1-\lambda)
		\left(
		abc+\frac{\gamma}{2}
		\right)+
		\lambda^2b^2\\
		&=
		(1-\lambda)^2a^2
		+2\lambda(1-\lambda)abc+
		\lambda(1-\lambda)\gamma
		+\lambda^2b^2.
	\end{align*}
    Together with (\ref{olambda2}), we have that
	\[
		4\rho^2\omega_\lambda^TC^{-1}\omega_\lambda
		=
		|\omega_\lambda|^2
		+\lambda(1-\lambda)\gamma = |\omega_\lambda|^2+q_\lambda\gamma. 
	\]
	Consequently, 
	\begin{align*}
	\ell(\lambda)^2
	&=
	|\omega_\lambda|^2+q_\lambda\gamma
	+
	\min_{|V|\leq1}
	\left(
	q_\lambda V^TCV
	+
	4\rho\sqrt{q_\lambda}\,
	\omega_\lambda^TV
	\right)\\
	&=
	4\rho^2\omega_\lambda^TC^{-1}\omega_\lambda + \min_{|V|\leq1}
	\left(
	q_\lambda V^TCV
	+
	4\rho\sqrt{q_\lambda}\,
	\omega_\lambda^TV
	\right)\\
	&=
	\min_{|V|\leq1}
	\Big(
	q_\lambda V^TCV
	+
	4\rho\sqrt{q_\lambda}\,
	\omega_\lambda^TV+
	4\rho^2
	\omega_\lambda^TC^{-1}\omega_\lambda
	\Big).
	\end{align*}
	For each $V$, we calculate
	\begin{align*}
		&~\quad
		\left(
		\sqrt{q_\lambda}\,V
		+
		2\rho C^{-1}\omega_\lambda
		\right)^T
		C
		\left(
		\sqrt{q_\lambda}\,V
		+
		2\rho C^{-1}\omega_\lambda
		\right)\\
		&=
		\left(\sqrt{q_\lambda}\,V\right)^T
		C
		\left(\sqrt{q_\lambda}\,V\right)+
		\left(\sqrt{q_\lambda}\,V\right)^T
		C
		\left(2\rho C^{-1}\omega_\lambda\right)+
		\left(2\rho C^{-1}\omega_\lambda\right)^T
		C
		\left(\sqrt{q_\lambda}\,V\right)+
		\left(2\rho C^{-1}\omega_\lambda\right)^T
		C
		\left(2\rho C^{-1}\omega_\lambda\right)\\
		&
		=
		q_\lambda V^TCV
		+
		4\rho\sqrt{q_\lambda}\,
		\omega_\lambda^TV
		+
		4\rho^2
		\omega_\lambda^TC^{-1}\omega_\lambda.
	\end{align*}
	Thus
	\[
		\ell(\lambda)^2
		=
		\min_{|V|\leq1}
		\left\|
		\sqrt{q_\lambda}\,V
		+
		2\rho C^{-1}\omega_\lambda
		\right\|_C^2.
	\]
	Because $\|\cdot\|_C\geq0$,
	\( r\longmapsto r^2 \)
	is increasing on $[0,\infty)$, and a norm is symmetric, 
	\begin{align*}
		\ell(\lambda)
		&=
		\sqrt{
			\min_{|V|\leq1}
			\left\|
			\sqrt{q_\lambda}\,V
			+
			2\rho C^{-1}\omega_\lambda
			\right\|_C^2
		}\\
		&=
		\min_{|V|\leq1}
		\left\|
		\sqrt{q_\lambda}\,V
		+
		2\rho C^{-1}\omega_\lambda
		\right\|_C\\
		&= \min_{V\in D}
		\left\|
		\sqrt{q_\lambda}\,V-p(\lambda) 
		\right\|_C \\
		&= \min_{V\in D}\left\|
		p(\lambda)-\sqrt{q_\lambda}\,V
		\right\|_C\\
		&=\operatorname{dist}_C
		\left(
		p(\lambda),
		\sqrt{q_\lambda}\,D
		\right).
	\end{align*}
	Since
	\(
	q_\lambda=\lambda(1-\lambda),
	\)
	we conclude that
	\begin{equation*}
		\ell(\lambda)
		=
		\operatorname{dist}_C
		\left(
		p(\lambda),
		\sqrt{\lambda(1-\lambda)}\,D
		\right).
	\end{equation*}
\end{proof}

    Since \(C\) is independent of \(\lambda\) and
    \[
    \omega_\lambda=(1-\lambda)e_0+\lambda e_1
    \]
    is affine in \(\lambda\), the function
    \[
    p(\lambda)=-2\rho C^{-1}\omega_\lambda
    \]
    is also affine. We now prove the dual formula needed for the convexity
    argument.

\begin{lemma}
	\label{lem:disk-distance-duality}
	Let \(C\in\mathbb{R}^{2\times 2}\) be symmetric positive definite, and define
	\[
	\|z\|_{C}:=\sqrt{z^{T}Cz},
	\qquad
	\|y\|_{C^{-1}}:=\sqrt{y^{T}C^{-1}y}.
	\]
	Let
	\(
	D=\{z\in\mathbb{R}^{2}:|z|\leq 1\}
	\)
	be the closed Euclidean unit disk. Then, for every
	\(p\in\mathbb{R}^{2}\) and every \(r\geq 0\),
	\[
	\operatorname{dist}_{C}(p,rD)
	=
	\max_{\|y\|_{C^{-1}}\leq 1}
	\left(
	y^{T}p-r|y|
	\right).
	\]
\end{lemma}

\begin{proof}
	
	First observe that the norms \(\|\cdot\|_{C}\) and
	\(\|\cdot\|_{C^{-1}}\) are dual to each other. Indeed,
	for every \(v,y\in\mathbb{R}^{2}\),
	\[
	y^{T}v
	=
	(C^{-1/2}y)^{T}(C^{1/2}v).
	\]
	By the Euclidean Cauchy--Schwarz inequality,
	\[
	y^{T}v
	\leq
	|C^{-1/2}y|\,|C^{1/2}v|
	=
	\|y\|_{C^{-1}}\|v\|_{C}.
	\]
	If \(v\neq0\), equality is attained at
	\[
	y=\frac{Cv}{\|v\|_C},
	\]
	while the case \(v=0\) is immediate.
	Consequently,
	\[
	\|v\|_{C}
	=
	\max_{\|y\|_{C^{-1}}\leq 1}y^{T}v.
	\]

	Suppose first that \(p\in rD\). Then
	\(
	\operatorname{dist}_{C}(p,rD)=0.
	\)
	For every \(y\in\mathbb{R}^{2}\), the Euclidean
	Cauchy--Schwarz inequality gives
	\[
	y^{T}p\leq |y|\,|p|\leq r|y|.
	\]
	Hence
	\[
	y^{T}p-r|y|\leq 0.
	\]
	Since \(y=0\) is admissible, it follows that
	\[
	\max_{\|y\|_{C^{-1}}\leq 1}
	\left(y^{T}p-r|y|\right)=0.
	\]
	Thus the formula holds when \(p\in rD\).
	
	Now suppose that \(p\notin rD\). Since \(rD\) is compact,
	there exists \(z_{*}\in rD\) such that
	\[
	d:=\operatorname{dist}_{C}(p,rD)
	=\|p-z_{*}\|_{C}>0.
	\]
	The point \(z_{*}\) minimizes the function
	\[
	\Phi(z):=\frac12\|p-z\|_{C}^{2}
	=\frac12(p-z)^{T}C(p-z)
	\]
	over the convex set \(rD\).
	
	For any \(z\in rD\) and \(0\leq t\leq 1\), define
	\(
	z(t):=z_{*}+t(z-z_{*}).
	\)
	By convexity of \(rD\), one has \(z(t)\in rD\). Since
	\(z_{*}\) is a minimizer,
	\[
	\left.\frac{d}{dt}\Phi(z(t))\right|_{t=0^{+}}\geq 0.
	\]
	A direct differentiation gives
	\begin{align*}
		\frac{d}{dt}\Phi(z(t))
		&=
		-(z-z_{*})^{T}
		C
		\bigl(p-z_{*}-t(z-z_{*})\bigr).
	\end{align*}
	Therefore,
	\(
	-(z-z_{*})^{T}C(p-z_{*})\geq 0,
	\)
	or equivalently,
	\(
	\bigl(C(p-z_{*})\bigr)^{T}(z-z_{*})\leq 0
    \)
	for every \(z\in rD\). 
	
	Define
	\[
	y_{*}:=\frac{C(p-z_{*})}{\|p-z_{*}\|_{C}}
	=\frac{C(p-z_{*})}{d}.
	\]
	Then
	\[
	\begin{aligned}
		\|y_{*}\|_{C^{-1}}^{2}
		&=
		y_{*}^{T}C^{-1}y_{*}=
		\frac{
			(p-z_{*})^{T}CC^{-1}C(p-z_{*})
		}{d^{2}}=
		\frac{
			(p-z_{*})^{T}C(p-z_{*})
		}{d^{2}}
		=1.
	\end{aligned}
	\]
	Thus \(y_{*}\) is admissible in the dual maximization.
	
	The preceding optimality condition becomes
	\[
	y_{*}^{T}(z-z_{*})\leq 0
	\qquad
	\text{for every }z\in rD.
	\]
	Hence
	\(
	\max_{z\in rD}y_{*}^{T}z=y_{*}^{T}z_{*}.
	\)
	On the other hand, the support function of the Euclidean
	disk \(rD\) is
	\(
	\max_{z\in rD}y^{T}z=r|y|.
	\)
    Consequently,
	\(
	y_{*}^{T}z_{*}=r|y_{*}|.
	\)
	It follows that
	\[
	\begin{aligned}
		y_{*}^{T}p-r|y_{*}|
		=
		y_{*}^{T}p-y_{*}^{T}z_{*}
		=
		y_{*}^{T}(p-z_{*})
	    =
		\frac{
			(p-z_{*})^{T}C(p-z_{*})
		}{d}
		=d.
	\end{aligned}
	\]
	Therefore,
	\begin{align}\label{dual1}
		\max_{\|y\|_{C^{-1}}\leq 1}
		\left(y^{T}p-r|y|\right)
		\geq d.
	\end{align}

	Conversely, for every \(y\) satisfying
	\(\|y\|_{C^{-1}}\leq 1\), we have
	\[
	\begin{aligned}
		y^{T}p-r|y|
		&=
		y^{T}p-\max_{z\in rD}y^{T}z
		=
		\min_{z\in rD}y^{T}(p-z)
		\leq
		\min_{z\in rD}\|y\|_{C^{-1}}\|p-z\|_{C}
		\leq
		\min_{z\in rD}\|p-z\|_{C}
		=
		d.
	\end{aligned}
	\]
	Taking the maximum over all admissible \(y\) yields
	\begin{align}\label{dual2}
		\max_{\|y\|_{C^{-1}}\leq 1}
		\left(y^{T}p-r|y|\right)
		\leq d.
	\end{align}
	Combining (\ref{dual1}) and (\ref{dual2}), we have that
 	\[
	\operatorname{dist}_{C}(p,rD)
	=
	\max_{\|y\|_{C^{-1}}\leq 1}
	\left(y^{T}p-r|y|\right).
	\qedhere
   \]
    
\end{proof}

We now apply this formula to the lower-boundary function.

\begin{proposition}
	\label{prop:lower-boundary-convexity}
	Suppose that the lower-boundary function has the representation
	\[
	\ell(\lambda)
	=
	\operatorname{dist}_{C}
	\left(
	p(\lambda),
	\sqrt{\lambda(1-\lambda)}\,D
	\right),
	\qquad 0\leq\lambda\leq 1,
	\]
	where \(C\) is symmetric positive definite,
	\(D=\{V\in\mathbb{R}^{2}:|V|\leq 1\}\), and
	\(p(\lambda)\) is affine in \(\lambda\). Then
	\(\ell\) is convex on \([0,1]\).
\end{proposition}

\begin{proof}
	Applying Lemma~\ref{lem:disk-distance-duality} with
	\(
	r(\lambda)=\sqrt{\lambda(1-\lambda)},
	\)
	we obtain
	\[
	\ell(\lambda)
	=
	\max_{\|y\|_{C^{-1}}\leq 1}
	\left[
	y^{T}p(\lambda)
	-
	|y|\sqrt{\lambda(1-\lambda)}
	\right].
	\]
	
	Fix \(y\) satisfying \(\|y\|_{C^{-1}}\leq 1\).
	Since \(p(\lambda)\) is affine in \(\lambda\), the function
	\(
	\lambda\longmapsto y^{T}p(\lambda)
	\)
	is affine.
	Moreover, on \(0<\lambda<1\),
	\(
	r(\lambda)=\sqrt{\lambda(1-\lambda)}
	\)
	satisfies
	\[
	r'(\lambda)
	=
	\frac{1-2\lambda}
	{2\sqrt{\lambda(1-\lambda)}}
	\]
	and
	\[
	\begin{aligned}
		r''(\lambda)
		=
		-\frac{1}{\sqrt{\lambda(1-\lambda)}}
		-
		\frac{(1-2\lambda)^{2}}
		{4[\lambda(1-\lambda)]^{3/2}}
	    =
		-\frac{
			4\lambda(1-\lambda)+(1-2\lambda)^{2}
		}{
			4[\lambda(1-\lambda)]^{3/2}
		}
	    =
		-\frac{1}
		{4[\lambda(1-\lambda)]^{3/2}}
		<0.
	\end{aligned}
	\]
	Thus \(r\) is concave on \([0,1]\). Consequently,
	\(
	\lambda
	\longmapsto
	-|y|\sqrt{\lambda(1-\lambda)}
	\)
	is convex. It follows that, for every fixed admissible \(y\), the function
	\[
	\lambda
	\longmapsto
	y^{T}p(\lambda)
	-
	|y|\sqrt{\lambda(1-\lambda)}
	\]
	is convex. Since the pointwise supremum of an arbitrary family
	of convex functions is convex, we conclude that
	\[
	\ell(\lambda)
	=
	\sup_{\|y\|_{C^{-1}}\leq 1}
	\left[
	y^{T}p(\lambda)
	-
	|y|\sqrt{\lambda(1-\lambda)}
	\right]
	\]
	is convex on \([0,1]\).	Finally, the real coordinate \(t\) is related to \(\lambda\) by
	\[
	t=(1-\lambda)h_{1}+\lambda h_{2}.
	\]
	If \(h_{1}<h_{2}\), then
	\[
	\lambda=\frac{t-h_{1}}{h_{2}-h_{1}}
	\]
	is an affine function of \(t\). Therefore the lower-boundary
	function expressed in the real coordinate \(t\) is also convex.
\end{proof}

\begin{remark}[Degenerate cases]
	\label{rem:degenerate-cases}
	The preceding calculation was carried out under the generic
	assumptions
	\[
	h_{1}<h_{2},
	\qquad
	\rho>0,
	\qquad
	a>0,
	\qquad
	b>0,
	\qquad
	s> 0.
	\]
	The remaining cases can be treated as follows.
	
	If \(h_{1}=h_{2}\), then the real coordinate of the upper bild
	takes only one value. Hence there is no nontrivial convexity
	statement to prove.
	
	If \(\rho=0\), then \(H\) and \(K\) are simultaneously block
	diagonal. On the fiber parametrized by \(\lambda\), the minimum
	is obtained by placing the two imaginary diagonal contributions
	in opposite directions. Consequently,
	\[
	\ell(\lambda)
	=
	\left|(1-\lambda)a-\lambda b\right|.
	\]
	This is the absolute value of an affine function and is therefore
	convex.
	
	Finally, suppose that one or more of \(a\), \(b\), or \(s\)
	vanishes. Choose a sequence of generic matrices \(K_{n}\)
	converging to \(K\). Such a sequence can be obtained by perturbing the vanishing
	normal-form parameters to positive values while keeping \(H\)
	fixed. Hence the fibers \(\mathcal F_\lambda\) are unchanged. For every \(\lambda\), write
	\[
	\ell_{K}(\lambda)
	:=
	\min_{x\in\mathcal{F}_{\lambda}}|x^{*}Kx|,
	\]
	where \(\mathcal{F}_{\lambda}\) denotes the corresponding compact
	fiber of quaternionic unit vectors.
	
	For every unit vector \(x\),
	\[
	\begin{aligned}
		\bigl|
		|x^{*}K_{n}x|-|x^{*}Kx|
		\bigr|
		\leq
		|x^{*}(K_{n}-K)x|
		\leq
		\|K_{n}-K\|.
	\end{aligned}
	\]
	Taking minima over \(x\in\mathcal{F}_{\lambda}\) and then
	interchanging \(K\) and \(K_{n}\), we obtain
	\[
	|\ell_{K_{n}}(\lambda)-\ell_{K}(\lambda)|
	\leq
	\|K_{n}-K\|.
	\]
	The right-hand side is independent of \(\lambda\); therefore
	\(\ell_{K_{n}}\to\ell_{K}\) uniformly on \([0,1]\).
	
	Each \(\ell_{K_{n}}\) is convex by the generic argument. Hence,
	for \(0\leq\theta\leq1\) and
	\(\lambda_{0},\lambda_{1}\in[0,1]\),
	\[
	\begin{aligned}
		\ell_{K_{n}}
		\bigl((1-\theta)\lambda_{0}+\theta\lambda_{1}\bigr)
		\leq
		(1-\theta)\ell_{K_{n}}(\lambda_{0})
		+
		\theta\ell_{K_{n}}(\lambda_{1}).
	\end{aligned}
	\]
	Passing to the limit \(n\to\infty\) gives
	\[
	\begin{aligned}
		\ell_{K}
		\bigl((1-\theta)\lambda_{0}+\theta\lambda_{1}\bigr)
		\leq
		(1-\theta)\ell_{K}(\lambda_{0})
		+
		\theta\ell_{K}(\lambda_{1}).
	\end{aligned}
	\]
	Thus the convexity conclusion also holds in the degenerate cases.
\end{remark}

\begin{theorem}\label{thm:main}
	For every quaternionic matrix \(A\in M_n(\mathbb H)\), the upper bild
	\(\bild^+(A)\) is convex.
\end{theorem}

\begin{proof}
	For \(n=1\), \(\bild^+(A)\) is a singleton, as observed in the Introduction, that is Section \ref{Section: Introduction}. We may therefore assume \(n\geq2\).
	For \(2\times2\) matrices, the upper-boundary function is concave by
	Lemma~\ref{upconcave}, while the lower-boundary function is convex by
	Proposition~\ref{prop:lower-boundary-convexity} and
	Remark~\ref{rem:degenerate-cases}. Together with
	Lemma~\ref{connected}, this proves the convexity of the upper bild in
	the \(2\times2\) case. The general case follows from
	Lemma~\ref{2x2}.
\end{proof}

\section{\textbf{Conclusion}}\label{Section:Conclusion}

In this paper, we studied the convexity and geometry of the upper bild of the
numerical range of a quaternionic matrix. For several special classes of \(2\times2\) quaternionic matrices,
we used the Schur form to obtain explicit geometric descriptions
of the shape and location of the upper bild in terms of the
standard eigenvalues and the off-diagonal parameters.
As another application of our analysis of the special cases,
we obtained some conditions for the quaternionic
numerical ranges of matrices in these special classes to be convex.

For the general \(2\times 2\) case, we employed a different canonical form and
developed an optimization-based approach to establish the convexity of the
upper bild. For the upper boundary, we related the corresponding maximization
problem to the complex numerical range of the associated complex
representation.
We proved that the upper boundary of the upper bild coincides
with that of the complex numerical range of the associated
complex representation, from which the concavity of the
upper-boundary function follows.
For the lower boundary, we obtained a distance-to-a-disk
representation and used it to prove the convexity of the
lower-boundary function. Combining the concavity of the upper-boundary function,
the convexity of the lower-boundary function, and the connectedness
of the vertical sections, we established the convexity of the
upper bild.

Several questions remain open. In particular, it would be interesting to
determine to what extent the algebraic and geometric properties of the upper
bild of a higher-dimensional quaternionic matrix can be recovered from those
of its \(2\times 2\) compressions. 

\section*{Acknowledgments}

The authors thank Professor Wasin So for drawing attention to gaps in Kumar's proof. In March 2021, he raised a question concerning Case 4 with the second author, who brought the problem to the first author's attention in January 2025. By the time of their in-person discussion with Professor So in October 2025, the authors had identified part of the gaps discussed in this paper. The authors are grateful to Professor So for his initial observation and subsequent discussions. We also thank Professor So for drawing our attention to the question of convexity for these special classes.

We became aware of the recent work of \cite{carvalho2026convexity}, posted on arXiv on October 13, 2026, after completing the first draft of this manuscript. Although both papers address the convexity of the upper bild of quaternionic numerical ranges, the present paper develops a different approach based on a distance-to-a-disk formulation and provides additional results on the boundary structure. Preliminary results on the gaps in Kumar's proof and the geometry of special classes were presented by the first author at ILAS 2026 in May 2026. The distance-to-a-disk approach was described in an abstract submitted for JMM 2027 in August 2026, before the appearance of Carvalho's preprint.

\section*{Declaration of generative AI and AI-assisted technologies
	in the manuscript preparation process}

During the preparation of this work, the authors used OpenAI's
ChatGPT to discuss mathematical arguments, organize the exposition,
and improve language and readability. They used OpenAI's Codex
and Anthropic's Claude Code to assist with exploratory computational
checks of selected formulas and examples. These checks were used
to look for possible inconsistencies, not as proofs of the stated
results. The authors reviewed and edited the AI-assisted output
and take full responsibility for the content of the publication.

\printcredits

\bibliographystyle{cas-model2-names}

\bibliography{Tian_LAA_2026_Upperbild}



\end{document}


\shorttitle{Supplementary material for the upper bild}
\shortauthors{J. Tian and T.-Y. Tam}

\title[mode=title]{Supplementary Mathematical Details for ``Revisiting the Upper Bild of Quaternionic Numerical Ranges''}

\author[1]{Jie Tian}[orcid=0009-0002-8793-2304]
\cormark[1]
\ead{jtian@unr.edu}
\ead[url]{https://jietinatian0803.github.io/}

\author[1]{Tin-Yau Tam}
\ead{ttam@unr.edu}
\ead[url]{https://sites.google.com/view/tinyautam/homepage}

\affiliation[1]{organization={University of Nevada, Reno},
  addressline={1664 N. Virginia Street},
  city={Reno},
  postcode={89557},
  state={Nevada},
  country={USA}}
\cortext[1]{Corresponding author}

\begin{abstract}
This supplement provides proofs of auxiliary reductions and computations
used in the main paper but omitted there to keep its main argument focused.
If requested by the editors or referees, the relevant proofs can be
incorporated into the main manuscript.
\end{abstract}

\maketitle

\appendix

\section{Canonical reductions}

We identify the space of pure quaternions with \(\R^3\).  Conjugation by a
unit quaternion induces a rotation of this space, and every rotation arises
in this way.

\begin{lemma}[Simultaneous normalization of two pure quaternions]
\label{supp:simultaneous-normalization}
Let \(\alpha,\beta\in\operatorname{Im}\H\).  There are
\(r\in\Sp(1)\), \(a,b\geq0\), \(c\in[-1,1]\), and \(s\geq0\) such that
\[
 r^*\alpha r=ai,
 \qquad
 r^*\beta r=b(ci+sj),
 \qquad c^2+s^2=1.
\]
If \(\alpha\neq0\) and \(\beta\neq0\), one may take
\[
 a=|\alpha|,
 \qquad b=|\beta|,
 \qquad
 c=\frac{\langle\alpha,\beta\rangle}{|\alpha|\,|\beta|},
 \qquad s=\sqrt{1-c^2}.
\]
\end{lemma}

\begin{proof}
Suppose first that \(\alpha\neq0\).  Choose \(r_0\in\Sp(1)\) such that
\(r_0^*\alpha r_0=|\alpha|i\), and write
\[
 r_0^*\beta r_0=\beta_i i+\beta_jj+\beta_kk.
\]
A rotation about the \(i\)-axis leaves \(i\) fixed and can rotate
\(\beta_jj+\beta_kk\) onto the nonnegative \(j\)-axis.  Hence there is
\(q\in\operatorname{span}_{\R}\{1,i\}\cap\Sp(1)\) for which
\[
 q^*iq=i,
 \qquad
 q^*(\beta_jj+\beta_kk)q
 =\sqrt{\beta_j^2+\beta_k^2}\,j.
\]
With \(r=r_0q\), the required representation follows.  Conjugation preserves
norms and Euclidean inner products, which gives the stated values of
\(a,b,c,s\).

If \(\alpha=0\) and \(\beta\neq0\), rotate \(\beta\) to \(|\beta|i\) and
take \(a=0\), \(b=|\beta|\), \(c=1\), and \(s=0\).  If \(\beta=0\), take
\(b=0\), \(c=1\), and \(s=0\); the remaining assertion follows by rotating
\(\alpha\) when it is nonzero.  The case \(\alpha=\beta=0\) is immediate.
\end{proof}

\begin{lemma}[Normalization of the Hermitian and anti-Hermitian parts]
\label{supp:matrix-normalization}
Let
\[
 H=\begin{pmatrix}h_1&0\\0&h_2\end{pmatrix},
 \qquad
 K=\begin{pmatrix}\alpha&z\\-z^*&\beta\end{pmatrix},
\]
where \(h_1,h_2\in\R\), \(\alpha,\beta\in\operatorname{Im}\H\), and
\(z\neq0\).  There is a diagonal quaternionic unitary matrix
\(D=\operatorname{diag}(r_1,r_2)\) such that
\[
 D^*HD=H,
 \qquad
 D^*KD=
 \begin{pmatrix}
 ai&\rho\\
 -\rho&b(ci+sj)
 \end{pmatrix},
\]
where \(\rho=|z|>0\), \(a=|\alpha|\), \(b=|\beta|\),
\(c\in[-1,1]\), \(s\geq0\), and \(c^2+s^2=1\).
\end{lemma}

\begin{proof}
Set \(\rho=|z|\) and
\[
 \widehat\beta=\frac{z\beta z^*}{\rho^2}.
\]
Then \(\widehat\beta\) is pure and \(|\widehat\beta|=|\beta|\).  By
Lemma~\ref{supp:simultaneous-normalization}, choose \(r_1\in\Sp(1)\) so that
\[
 r_1^*\alpha r_1=ai,
 \qquad
 r_1^*\widehat\beta r_1=b(ci+sj).
\]
Define \(r_2=z^*r_1/\rho\).  Then \(r_2\in\Sp(1)\) and
\[
 r_1^*zr_2=\rho,
 \qquad
 r_2^*\beta r_2
 =r_1^*\frac{z\beta z^*}{\rho^2}r_1
 =b(ci+sj).
\]
This proves the asserted form of \(D^*KD\).  Since the entries of \(H\) are
real, every diagonal unitary matrix commutes with \(H\), and hence
\(D^*HD=H\).
\end{proof}

\begin{lemma}[Normalization of the off-diagonal entry]
\label{supp:triangular-normalization}
Let
\[
 A=\begin{pmatrix}
 a_1+b_1i&q\\0&a_2+b_2i
 \end{pmatrix},
 \qquad b_1,b_2\geq0,
\]
and write \(q=z_1+z_2j\), where \(z_1,z_2\in\C\).  Then \(A\) is
unitarily similar to a matrix of the form
\[
 \begin{pmatrix}
 a_1+b_1i&m+\gamma j\\0&a_2+b_2i
 \end{pmatrix},
 \qquad m=|z_1|,\quad \gamma=|z_2|.
\]
\end{lemma}

\begin{proof}
Choose complex units \(r_1,r_2\in\C\cap\Sp(1)\).  Because they commute with
\(i\), the diagonal unitary \(D=\operatorname{diag}(r_1,r_2)\) leaves both
diagonal entries unchanged.  Its action on the off-diagonal entry is
\[
 r_1^*qr_2
 =\overline r_1z_1r_2+\overline r_1z_2\overline r_2j.
\]
Write \(z_1=me^{i\varphi}\) and \(z_2=\gamma e^{i\psi}\), with arbitrary
phases when one of the two numbers is zero.  It is enough to choose phases
\(r_1=e^{i\theta_1}\), \(r_2=e^{i\theta_2}\) satisfying
\[
 -\theta_1+\varphi+\theta_2=0,
 \qquad
 -\theta_1+\psi-\theta_2=0.
\]
This real linear system always has a solution.  For such a choice,
\(r_1^*qr_2=m+\gamma j\), as required.
\end{proof}

\begin{lemma}[Relation between the two generic cases]
\label{supp:generic-cases}
Let
\[
 A=
 \begin{pmatrix}
 a_1+b_1i&m+\gamma j\\
 0&a_2+b_2i
 \end{pmatrix},
 \qquad b_1,b_2>0,
 \qquad m,\gamma\geq0,
 \qquad m+\gamma>0.
\]
Write \(A=H+K\), where \(H^*=H\) and \(K^*=-K\).  After diagonalizing
\(H\), put \(K\) in the normal form of Lemma~\ref{supp:matrix-normalization}
\[
 \begin{pmatrix}
 ai&\rho\\
 -\rho&b(ci+sj)
 \end{pmatrix},
 \qquad a,b,\rho\geq0,
 \qquad s\geq0,
 \qquad c^2+s^2=1.
\]
Then
\[
 m>0\ \text{and}\ \gamma>0
 \quad\Longleftrightarrow\quad
 a>0,\ b>0,\ \rho>0,\ \text{and}\ s>0.
\]
Thus the remaining case \(m>0,\gamma>0\) in the triangular form is
exactly the generic case used in the Hermitian--anti-Hermitian reduction.
\end{lemma}

\begin{proof}
Set
\[
 q=m+\gamma j,
 \qquad
 r=|q|=\sqrt{m^2+\gamma^2},
 \qquad
 w=\frac{q}{r}.
\]
Since
\[
 H=
 \begin{pmatrix}
 a_1&q/2\\
 q^*/2&a_2
 \end{pmatrix},
 \qquad
 K=
 \begin{pmatrix}
 b_1i&q/2\\
 -q^*/2&b_2i
 \end{pmatrix},
\]
the diagonal unitary matrix \(P=\operatorname{diag}(w,1)\) gives
\[
 P^*HP=
 \begin{pmatrix}
 a_1&r/2\\
 r/2&a_2
 \end{pmatrix},
 \qquad
 P^*KP=
 \begin{pmatrix}
 b_1f&r/2\\
 -r/2&b_2i
 \end{pmatrix},
\]
where
\[
 f=w^*iw
 =u i+v k,
 \qquad
 u=\frac{m^2-\gamma^2}{r^2},
 \qquad
 v=\frac{2m\gamma}{r^2}.
\]

Choose a real orthogonal matrix
\[
 R=
 \begin{pmatrix}
 C&-S\\
 S&C
 \end{pmatrix},
 \qquad C^2+S^2=1,
\]
that diagonalizes \(P^*HP\).  Since \(r>0\), neither \(C\) nor \(S\)
is zero.  A direct multiplication gives
\[
 R^T(P^*KP)R
 =
 \begin{pmatrix}
 \alpha&z\\
 -z^*&\beta
 \end{pmatrix},
\]
where
\[
 \alpha=C^2b_1f+S^2b_2i,
 \qquad
 \beta=S^2b_1f+C^2b_2i,
 \qquad
 z=\frac r2+CS(b_2i-b_1f).
\]
In particular, \(\operatorname{Re}z=r/2>0\), so \(\rho=|z|>0\).

Suppose first that \(m>0\) and \(\gamma>0\).  Then \(v\neq0\).  The
\(k\)-components of \(\alpha\) and \(\beta\) are \(C^2b_1v\) and
\(S^2b_1v\), respectively.  Hence \(\alpha\neq0\) and \(\beta\neq0\),
so \(a=|\alpha|>0\) and \(b=|\beta|>0\).

It remains to verify that the two pure quaternions used in the simultaneous
normalization are not parallel.  As in
Lemma~\ref{supp:matrix-normalization}, these quaternions are
\[
 \alpha
 \quad\text{and}\quad
 \widehat\beta=\frac{z\beta z^*}{|z|^2}.
\]
Identify pure quaternions with vectors in \(\R^3\).  Direct calculation
from the preceding formulas shows that the \(i\)-component of
\(\alpha\times(z\beta z^*)\) is
\[
 b_1^2b_2 C^3Sr v^2.
\]
This number is nonzero because \(b_1,b_2,r>0\), \(CS\neq0\), and
\(v\neq0\).  Thus \(\alpha\) and \(\widehat\beta\) are not parallel,
and the parameter \(s\) in the normal form is positive.

Conversely, suppose that \(m\gamma=0\).  Then \(v=0\) and \(f=\pm i\).
Consequently, \(\alpha\) and \(\beta\) are real multiples of \(i\), while
\(z\in\operatorname{span}_{\R}\{1,i\}\).  Conjugation by \(z/|z|\)
therefore leaves the real line spanned by \(i\) invariant.  Hence
\(\alpha\) and \(z\beta z^*/|z|^2\) are parallel whenever both are
nonzero; if one of them is zero, then \(a=0\) or \(b=0\).  In either
case the generic conditions \(a,b,\rho,s>0\) fail.  This proves the
equivalence.
\end{proof}

\begin{remark}[Interchanging a real and a nonreal diagonal entry]
\label{supp:order}
The case in which the nonreal standard eigenvalue appears first and the real
eigenvalue second does not give a new upper bild.  Indeed,
\(W_{\H}(A^*)=W_{\H}(A)\).  After taking the adjoint, a permutation unitary
interchanges the two diagonal positions.  Finally, quaternionic conjugation
within a similarity class replaces a complex diagonal entry by its standard
representative in \(\C^+\).  Thus the order used in the main paper entails no
loss of generality.
\end{remark}

\section{A further remark on \cite{K19}}
\label{supp:kumar}

In this appendix we discuss a statement in \cite{K19} that is not
correct as stated. This does not affect the convexity of the upper bild,
and the statement is not used in the main paper. We use the notation of \cite{K19}:
\(\langle X,Y\rangle_{\H}=\sum_{l}\overline{x_l}\,y_l\),
\(\mathcal S_{\H^n}\) is the unit sphere of \(\H^n\), and
\(\operatorname{co}(q)=q_0+q_1i\) for
\(q=q_0+q_1i+q_2j+q_3k\).

\subsection{On Theorem~2.7 of \cite{K19}}

Theorem~2.7 of \cite{K19} asserts that, for \(A\in M_n(\H)\), the
convexity of \(W_{\H}(A)\) is equivalent to the following condition:
\begin{itemize}
\item[(4)] for every \(X\in\mathcal S_{\H^n}\) there exists
\(Y\in\mathcal S_{\H^n}\) such that
\[
2\,\bigl|\operatorname{im}\langle Y,AY\rangle_{\H}\bigr|
=\bigl|\langle X,AX\rangle_{\H}\,i-i\,\overline{\langle X,AX\rangle_{\H}}\bigr| .
\]
\end{itemize}

\begin{lemma}\label{supp:co-identity}
For every \(q=q_0+q_1i+q_2j+q_3k\in\H\),
\(
q\,i-i\,\overline q=-2q_1 .
\)
In particular, the right-hand side of~(4) equals
\(2\,|\operatorname{im}\operatorname{co}\langle X,AX\rangle_{\H}|\).
\end{lemma}

\begin{proof}
We have \(q\,i=-q_1+q_0i+q_3j-q_2k\) and
\(i\,\overline q=q_1+q_0i+q_3j-q_2k\). Subtracting gives the claim.
\end{proof}

\begin{proposition}\label{supp:thm27}
The implication \((4)\Rightarrow(1)\) in \cite[Theorem~2.7]{K19} is false.
More precisely, for
\[
A=\begin{pmatrix} i&0\\0&1\end{pmatrix},
\]
condition~(4) holds, but \(W_{\H}(A)\) is not convex.
\end{proposition}

\begin{proof}
\emph{\(W_{\H}(A)\) is not convex.} Taking \(X=(1,0)^T\) and
\(X=(j,0)^T\) gives \(i\in W_{\H}(A)\) and
\(\overline{j}\,i\,j=-i\in W_{\H}(A)\). If \(0\in W_{\H}(A)\), then
\(0=\overline{x_1}\,i\,x_1+|x_2|^2\) for some unit vector
\((x_1,x_2)^T\). The first term is purely imaginary with modulus
\(|x_1|^2\), and the second is real, so \(x_1=0\) and \(x_2=0\), which is
impossible. Hence \(0=\frac12 i+\frac12(-i)\notin W_{\H}(A)\). This is
also Example~3.4 of~\cite{So2023}.

\emph{Condition~(4) holds.} Let \(X=(x_1,x_2)^T\in\mathcal S_{\H^2}\) and
\(q=\langle X,AX\rangle_{\H}=\overline{x_1}\,i\,x_1+|x_2|^2\). Then
\(|q_1|\le|\operatorname{im}q|=|x_1|^2\le1\). Choose
\(\sigma\in[0,\frac{\pi}{2}]\) with \(\cos^2\sigma=|q_1|\), and put
\(Y=(\cos\sigma,\sin\sigma)^T\). Then
\(\langle Y,AY\rangle_{\H}=\cos^2\sigma\,i+\sin^2\sigma\), so
\[
2\,|\operatorname{im}\langle Y,AY\rangle_{\H}|=2\cos^2\sigma=2|q_1|
=|q\,i-i\,\overline q|
\]
by Lemma~\ref{supp:co-identity}.
\end{proof}

\begin{remark}
The proof of \((4)\Rightarrow(2)\) in \cite{K19} uses the equality
\(\operatorname{re}\langle Y,AY\rangle_{\H}=\operatorname{re}\langle X,AX\rangle_{\H}\),
which is obtained in the proof of \((3)\Rightarrow(4)\) but is not part
of condition~(4). If this equality is added to~(4), the resulting
condition is equivalent to~(2), by the argument given in \cite{K19}.
\end{remark}

\clearpage
\section{Suggested parameters for the interactive visualization}
\label{supp:visualization-parameters}

Table~\ref{tab:supp-visualization-parameters} gives suggested inputs for
the accompanying file \texttt{upper\_bild\_interactive.html}, organized
by the sections of the main paper. The matrix entered in the interface is
\[
A=\begin{pmatrix}
a_1+b_1i&q_0+q_1i+q_2j+q_3k\\
0&a_2+b_2i
\end{pmatrix}.
\]
For every row, set \(q_1=q_3=0\). The values of \(q_0\) and \(q_2\)
are given in the table; in Section~4.4 these are \(m\) and \(\gamma\),
respectively. For reproducible plots, select 200000 samples and random
seed 2026, then click \texttt{Compute upper bild}. Finite sampling
illustrates the geometric descriptions proved in the main paper; it
does not determine exact boundary points.

\begin{table}[pos=!h]
\centering
\caption{Suggested parameter sets for Sections~4.1--4.4 of the main paper.
All unspecified quaternionic components are zero.}
\label{tab:supp-visualization-parameters}
\small
\setlength{\tabcolsep}{4pt}
\renewcommand{\arraystretch}{1.18}
\begin{tabular}{@{}lccc>{\raggedright\arraybackslash}p{0.45\textwidth}@{}}
\toprule
\shortstack{Main\\section} & \((a_1,b_1,a_2,b_2)\) & \(q_0\) & \(q_2\)
& Upper-bild geometry \\
\midrule
4.1 & \((0,1,2,1)\) & 0 & 0
& Triangle with vertices \(i\), \(2+i\), and \(1\). \\
4.1 & \((0,1,0,2)\) & 0 & 0
& Vertical segment from \(0\) to \(2i\). \\
4.1 & \((0,0,2,0)\) & 0 & 0
& Real segment \([0,2]\). \\
\hline
4.2 & \((-1,0,1,0)\) & 2 & 0
& Upper half of an elliptical disk, with horizontal semiaxis
\(\sqrt{2}\) and vertical semiaxis \(1\). \\
4.2 & \((0,0,0,0)\) & 2 & 0
& Upper half of the unit disk. \\
\hline
4.3 & \((0,0,2,1)\) & 2 & 0
& Tilted elliptical disk with foci \(0\) and \(2+i\), restricted
to the closed upper half-plane. \\
\hline
4.4.1 & \((0,1,2,1)\) & 1 & 0
& \(m^2<4b_1b_2\): the lower boundary includes two line segments. \\
4.4.1 & \((0,1,2,1)\) & 2 & 0
& \(m^2=4b_1b_2\): the ellipse is tangent to the real axis. \\
4.4.1 & \((0,1,2,1)\) & 3 & 0
& \(m^2>4b_1b_2\): the elliptical disk is cut by the real axis. \\
\hline
4.4.2 & \((0,1,2,4)\) & 0 & 4
& \(\Lambda=-80<0\): no additional upper common tangent segment. \\
4.4.2 & \((0,1,3,4)\) & 0 & 4
& \(\Lambda=0\): the critical case, with no nontrivial upper
common tangent segment. \\
4.4.2 & \((0,1,4,4)\) & 0 & 4
& \(\Lambda=112>0\): the upper boundary includes a common
tangent segment. \\
\bottomrule
\end{tabular}
\end{table}
\clearpage

\bibliographystyle{cas-model2-names}
\bibliography{Tian_LAA_2026_Upperbild}